\documentclass[seceqn,MSNbibl,number,citesort,dvips]{arxbj}
\usepackage{mathbh}

\aid{0}
\volume{20}
\issue{3}
\pubyear{2014}
\firstpage{1560}
\lastpage{1599}
\doi{10.3150/13-BEJ533} 

\makeatletter

\newcommand{\rright}{\right}
\newcommand{\lleft}{\left}
\newcommand{\rrVert}{\Vert}
\newcommand{\rrvert}{\vert}
\newcommand{\llVert}{\Vert}
\newcommand{\llvert}{\vert}

\newcommand{\cA}{\mathcal{A}}
\newcommand{\cC}{\mathcal{C}}
\newcommand{\cD}{\mathcal{D}}
\newcommand{\cF}{\mathcal{F}}
\newcommand{\cH}{\mathcal{H}}
\newcommand{\cK}{\mathcal{K}}
\newcommand{\cL}{\mathcal{L}}
\newcommand{\cM}{\mathcal{M}}
\newcommand{\cN}{\mathcal{N}}
\newcommand{\cP}{\mathcal{P}}
\newcommand{\cS}{\mathcal{S}}
\newcommand{\cV}{\mathcal{V}}
\newcommand{\cY}{\mathcal{Y}}
\newcommand{\bG}{\mathbb G}
\newcommand{\bH}{\mathbb H}
\newcommand{\bN}{{\mathbb N}}
\newcommand{\bP}{{\mathbb P}}
\newcommand{\bR}{{\mathbb R}}
\newcommand{\mr}{\mathbb{R}}
\newcommand{\me}{\mathbb{E}}

\newcommand{\op}[1]{\operatorname{#1}} 

\newtheorem{theorem}{Theorem}

\newtheorem{proposition}{Proposition}
\newtheorem{lemma}{Lemma}
\newproclaim{definition}{Definition}
\newremark{remark}{Remark}
\newremark{Example}{Example}
\newtheorem{corollary}{Corollary}
\newproclaim{condition}{Condition}
\newproclaim{discussion}{Discussion}

\newcommand{\1}{\mathbh{1}}

\newcommand{\penn}{\op{ani}_\epsilon(n)}

\newcommand{\penni}{\op{iso}_\epsilon(n)}

\newcommand{\netiso}{\cH_\epsilon^{\op{iso}}}

\newcommand{\eigen}{\Phi_{h}}

\newcommand{\estimiso}{\hat f_{\op{iso}}}

\newcommand{\eqref}[1]{(\ref{#1})}

\newcommand{\firstestimator}[1]{\hat{f}_{#1}}
\newcommand{\xrightarrow}[1]{\stackrel{#1}{\rule[2.3pt]{10pt}{0.3pt}\hspace*{-4pt}\longrightarrow}}

\def\sfrac#1#2{#1/#2}
\def\vfrac#1#2{(#1)/#2}

\def\sklfrac#1#2{(#1/#2)}
\makeatother

\begin{document}
\begin{frontmatter}

\title{A robust, adaptive M-estimator for pointwise estimation in
heteroscedastic regression}
\runtitle{A fully adaptive pointwise M-estimator}

\begin{aug}
\author{\inits{M.}\fnms{Micha\"el} \snm{Chichignoud}\corref{}\thanksref{e1}\ead[label=e1,mark]{chichignoud@stat.math.ethz.ch}} \and
\author{\inits{J.}\fnms{Johannes} \snm{Lederer}\thanksref{e2}\ead[label=e2,mark]{lederer@stat.math.ethz.ch}}
\runauthor{M. Chichignoud and J. Lederer} 
\address{Seminar for Statistics, ETH Z\"urich, R\"amistrasse 101,
CH-8092 Z\"urich, Switzerland.\\ \printead{e1,e2}}
\end{aug}

\received{\smonth{1} \syear{2013}}
\revised{\smonth{5} \syear{2013}}

%
\begin{abstract}
We introduce a robust and fully adaptive method for pointwise estimation
in heteroscedastic regression. We allow for noise and design distributions
that are unknown and fulfill very weak assumptions only. In particular, we
do not impose moment conditions on the noise distribution. Moreover, we do
not require a positive density for the design distribution. In a first
step, we study the consistency of locally polynomial M-estimators that
consist of a
contrast and a kernel. Afterwards, minimax results are established over
unidimensional H\"{o}lder spaces for degenerate design. We then choose
the contrast and the kernel that
minimize an empirical variance term and demonstrate that the
corresponding M-estimator is adaptive with respect to the noise and
design distributions and adaptive (Huber) minimax for contamination
models. In a second step, we additionally choose a data-driven bandwidth
via Lepski's method. This leads to an M-estimator that is adaptive with
respect to the noise and design distributions and, additionally,
adaptive with respect to the smoothness of an isotropic, multivariate, locally
polynomial target function. These results are also extended to
anisotropic, locally constant target functions. Our data-driven
approach provides, in particular, a level of robustness that adapts to
the noise, contamination, and outliers.
\end{abstract}

%

\end{frontmatter}
%
\section{Introduction}\label{intro}
We introduce a new method for pointwise estimation in heteroscedastic
regression that is adaptive with respect to the model, in particular,
with respect to the noise and the design distribution (D-adaptive) and
the smoothness of the
regression function (S-adaptive).

Let us first briefly summarize the related literature. First, the
seminal paper \cite{Huber64} contains a proof of the asymptotic normality
of M-estimators for the location parameter in regular models.
Furthermore, the series of papers \cite
{Tsybakov82a,Tsybakov82b,Tsybakov83,Tsybakov86} provide minimax results for
nonparametric regression. More recently, a
block median method was used in \cite{Cai_Zhou09} to prove the
asymptotic equivalence between Gaussian regression and homoscedastic
regression for deterministic designs and possibly heavy-tailed
noises. Using a blockwise Stein's Method with wavelets, this leads to
an S-adaptive estimator that is adaptive optimal over Besov spaces with
respect to the $ L_2 $-risk and adaptive optimal over isotropic H\"{o}lder
classes with respect to the punctual risk. Moreover, using an estimate
of the noise density at $0$ and a plug-in method, this also leads to a
D-adaptive estimator. However,
in contrast to this paper, only homoscedastic regression is considered
and multivariate regression functions, in particular anisotropic
functions, are not allowed for. Next, a modified version of Lepski's
method was applied for homoscedastic regression in \cite
{Reiss_Rozenholc_Cuenod09}. Finally, local M-estimators, also for
regression models with degenerate designs, were intensively studied in
the case of Gaussian regression: S-adaptivity results of a local least
squares estimator were derived in \cite{Gaiffas07}, sup-norm S-minimax
results were established in \cite{Gaiffas07b}, and the effect of
degenerate designs on the $ L_2 $-norm was investigated with
wavelet-type estimators in \cite{Antoniadis_Pensky_Sapatinas13}.
However, in contrast to this paper, pointwise estimation with random,
possibly degenerate designs and heteroscedastic, possibly heavy-tailed
noises has not been included.
%

What is the main idea behind our approach? Consider the estimation of $
t^0\in\bR$ in the translation model $ \cY\sim
g(\cdot-t^0) $ for a probability density $g$. The M-estimator $ \hat t
$ of $ t^0 $ corresponding to the contrast $\rho(\cdot)$ and the sample
$ \cY_1,\ldots,\cY_n $ of $ \cY$ is then
\[
\hat t:=\arg\min_t \sum_{i=1}^n
\rho(\cY_i-t).
\]
It holds that (see \cite{Huber64,Huber81,Huber_Ronchetti09})
%
\begin{equation}
\label{result normality huber} \sqrt n\bigl(\hat t-t^0\bigr)\mathop{\xrightarrow{\cL}}_{n\rightarrow
\infty}\cN(0, \mathrm{AV}), \qquad \mbox{where }\mathrm{AV}:=
\frac{\int(\rho
')^2\,\mathrm{d}G}{ (\int\rho''\,\mathrm{d}G )^2} ,
\end{equation}
$ G $ is the distribution of $\cY-t^0$, $ \rho'(\cdot) $ and $
\rho''(\cdot) $ are the first and second derivatives of the contrast
$\rho(\cdot)$, and $ \cL$ indicates convergence in law. In other
words, $ \hat t $ is asymptotically normal with asymptotic variance
$\mathrm{AV}$. This result suggests that an optimal estimator is obtained
by minimizing the asymptotic variance. Moreover, the {Cr\'{a}mer--Rao
Inequality} and (see \cite{Huber64})
%
\begin{equation}
\label{result fisher huber} \inf_\rho \frac{\int(\rho')^2\,\mathrm{d}G}{ (\int\rho''\,\mathrm{d}G )^2}= \bigl(I(G)
\bigr)^{-1},
\end{equation}
where $ I(\cdot) $ is the {Fisher information} and the
infimum is taken over all twice differentiable contrasts, imply that
this M-estimator is efficient. Huber proposed in \cite{Huber64}, Proposal
3, to minimize an estimate of the above asymptotic variance
(since the distribution $G$ is not available in practice) over the
family of Huber contrasts (their definition is given below). He also
conjectured that the corresponding estimator is minimax for certain
contamination models (for more details, see Section A.1 in the arXiv
version). More recently, in \cite{Arcones05}, an M-estimator with a
contrast that minimizes an estimate of the asymptotic variance was
introduced for the parametric model, its asymptotic normality was
proved, and especially Huber contrasts indexed by their scale and a
family of $\ell_p$ losses were considered.

In a first step, we derive general properties of M-estimators such as
pointwise risk bounds. This includes, in particular, S-minimax results
for degenerate designs and allows us to recover results in \cite{Gaiffas05}
(see Theorem~\ref{thr risk bound fix contrast} and Remark~\ref{remark
th 1 constant}). In a second step, we then consider a local M-estimator
that consists of a contrast and a kernel that minimize an estimate of
the variance and show, in particular, that this estimator mimics the
oracle, which minimizes the true variance. Our data-driven approach can
be used, for example, for the selection of the scale of the Huber
contrast with an
adaptive robustness with respect to outliers or for the selection of a
suitable (even noncentered or nonconvex) support that takes a maximal
number of
points around $x_0$ into account
(cf.\ \cite{Goldenshluger_Nemirovski97} for the latter objective).
Finally, we show that our
estimator is, under some restrictions on the design and the noise level (see
Condition \ref{condition variance estimation}), D-adaptive for various
sets of contrasts and kernels with
finite entropy.

We finally study simultaneous D- and S-adaptation for anisotropic
target functions. In a first step, we study the case of isotropic
target functions, where the standard
Lepski's method (see \cite{Lepski90,Lepski_Mammen_Spokoiny97}) can be
applied. To this end, we assume that the variance of the estimator is
decreasing with respect to the bandwidth and plug-in an estimate of the
minimal variance for the D-adaptation to apply Lepski's method for the
S-adaptation (see Section~\ref{section adaptation iso}). This yields
the first estimator in
heteroscedastic regression with random designs and heavy-tailed noise
distributions that is simultaneously D- and S-adaptive and optimal in a
sense describe later. Furthermore, we note that applications of
Lepski's method to nonlinear estimators are still nonstandard and
can only be found in a small number of examples in the literature
\cite{Chichignoud10,Polzehl_Spokoiny06,Reiss_Rozenholc_Cuenod09}. In a next step, we extend
our results to anisotropic target functions. For this, we restrict
ourselves to locally constant target functions and homoscedastic
regression with uniform design and apply a modification of
Lepski's method given in \cite
{Lepski_Levit99,Kerkyacharian_Lepski_Picard01} to construct an optimal,
simultaneously
S- and D-adaptive estimator. This is the first application of Lepski's method
to nonlinear estimators of anisotropic target functions and yields a
selection of an anisotropic bandwidth which is of great interest
for applications in the context of image denoising
(cf.\ \cite{Katkovnik_Foi_Egiazarian_Astola10}), for example.

Although we consider estimation problems, our approach may also be
useful for inference, for example, for the construction of confidence
bands. While confidence bands for parametric estimation are derived
from central limit theorems (see \eqref{result normality huber}),
confidence bands for nonparametric regression are especially desired to
be adaptive with respect to the smoothness of the target function. The
construction of such {S-adaptive confidence bands} is more difficult
than in the parametric case (see \cite{Hoffmann_Nickl11}), but since
Lepski-type procedures have already been used in this context, see
\cite{Gine_Nickl10}, Theorem~1 and Corollary~1, we expect that our
approach may be useful for the construction of S-adaptive confidence
bands for regression with possibly heavy-tailed noises (see Section~\ref{section discussion} for a discussion of some technical aspects).
Eventually, if for example the smoothness is known, our approach may be
used, plugging an estimate of the variance in the confidence band, to
obtain D-adaptive confidence bands, which are, in particular, adaptive
with respect to the design and the noise distributions.

The structure of this paper is as follows: In the following section, we
first introduce an estimator which satisfies a risk bound
(see Theorem~\ref{thr risk bound fix contrast}). So, S-minimax results
are deduced over H\"{o}lder spaces (see Corollaries \ref{coro
minimax}, \ref{coro minimax 2} and \ref{coro minimax 3}). We then
provide a choice
for the contrast and the kernel (see Theorem~\ref{thr oracle bound}) via
the minimization of a nonasymptotic variance. Then, we provide a
choice for the bandwidth for isotropic, locally polynomial target
functions (see Theorem~\ref{Th adaptive iso}) and for anisotropic,
locally constant target functions (see Theorem~\ref{Th adaptive aniso}).
After this, we give a discussion on our assumptions and an outlook in
Section~\ref{section discussion}. The proofs are finally conducted in
Section~\ref{section proofs main res} and in the \hyperref[app]{Appendix}. For
conciseness, only the crucial proofs are presented here. For the
remaining proofs and more details, in particular, on the parametric
model and on a comparison to classical results, we refer to the longer
version available on arXiv and the webpages of the authors.

\section{Preliminary definitions and results}

In this section, we give some preliminary definitions and results. After
specifying the model, we introduce a first estimator and then, we
present a risk bound and S-minimax properties of this
estimator.

Let us first specify the model. The observations
$(X_i,Y_i)_{i=1,\ldots,n}$ satisfy the set of equations
%
\begin{equation}
\label{model} Y_i=f^*(X_i)+\sigma(X_i)
\xi_i,\qquad  i=1,\ldots, n,
\end{equation}
and are distributed according to the probability measure $ \bP:=
\bP_{f^*}^{(n)}$ with associated\vspace*{-1pt} expectation $ \mathbb E:=\mathbb
E_{f^*}^{(n)} $.
We aim at estimating the target function
$f^*\dvtx [0,1]^d\rightarrow[ -M,M]$\vspace*{-1pt} (for $M>0$) at a given point $x_0$
on $ (0,1)^d $. The target function is assumed to be smooth, more
specifically, it is assumed to belong to a H\"older class (see
Definition~\ref{def holder space} below). The target function is
obscured by the second part of the above model, the noise. The noise
variables
$(\xi_i)_{i\in{1,\ldots,n}}$ are assumed to be distributed
independently according to the densities $g_i(\cdot) $ with respect to
the Lebesgue measure on $\mathbb R$. The noise densities $g_i(\cdot) $
may be unknown but are assumed to be symmetric. We stress that we do not
impose, unlike in the literature on the median (cf.\ \cite{Cai_Zhou09}),
any moment assumptions on the noise, and we do not require that the
noise densities are positive at $0$. We postpone the detailed discussion
on the assumptions to the end of the next section. The noise level $
\sigma\dvtx [0,1]^d\rightarrow[0,\infty) $ is assumed to be bounded, but may
also be unknown. Usually, the noise level is the variance of the noise,
however, this is not the case if the noise distributions do not have any
moments, for example. Finally, the design points $(X_i)_{i\in{1,\ldots,n}}$
are assumed to be distributed independently and identically according to
the density $\mu(\cdot)$ with respect to
the Lebesgue measure on $\mathbb R$. We assume that $\mu(\cdot)$
vanishes at most finitely many points. For ease of exposition, we
also assume that
$(X_i)_{i\in{1,\ldots,n}}$ and $(\xi_i)_{i\in{1,\ldots,n}}$ are mutually
independent.





Next, we introduce an estimator of $f^*(x_0)$ with a local polynomial
approach (LPA) for a fixed bandwidth, a fixed kernel, and a fixed contrast.
The key idea of the LPA, as described for example in \cite{Katkovnik85}
or in \cite{Tsybakov08}, Chapter~1, is to
approximate the target function in a neighborhood of size $
h\in(0,1]^d $ of a given point $x_0$ by a polynomial. To start, we
define for a fixed $m\in\bN$ the set $\cP:=\{
p=(p_1,\ldots,p_d)^\top\in\bN^d\dvtx  0\leq|p|
\leq m\}$ with $ |p|=p_1+\cdots+p_d $ and denote its cardinality by
$|\cP|$. The cardinality $|\cP|$ is exponential in $d$ and enters the
bounds derived below as a factor.
For any multi-indexed column vector ${t}= (t_{p_1,\ldots,p_d}\in
\bR\dvtx
p\in\cP )\in\bR^{|\cP|} $ and for any $ x\in[0,1]^d$, we then define
the desired polynomial as
\begin{eqnarray*}
\label{def polynomial} \mathrm{P}_t(x):= t^\top U \biggl(
\frac{x-{x_0}}{h} \biggr):=\sum_{p\in\cP}t_{p}
\biggl(\frac{x-{x_0}}{h} \biggr)^p.
\end{eqnarray*}
Here, $z^p:=z_1^{p_1}\cdots z_d^{p_d} $ for all $z\in\mr^d$, and the
division by $h$ is understood coordinate wise. Next, for $M>0$, we
define $
\cF:= \{\mathrm{P}_t \dvtx  t\in
[-M,M]^{|\cP|} \} $ as a set of polynomials of degree at most $
m $. We now specify what we mean by a kernel and a contrast:
%
\begin{definition}\label{def kernel}
A function $ K\dvtx \bR^d\rightarrow[0,\infty) $ is called \textit{kernel}
(function) if it has the
following properties:
\begin{enumerate}[2.]
\item$ K(\cdot)$ has a (not necessarily symmetric) support which is
a hypercube having edge length one and contains the origin;
\item$ \|K\|_\infty<\infty$ and $\int K(x)\,\mathrm{d}x=1 $.
\end{enumerate}
\end{definition}

For ease of exposition, we set $\Pi_h:=\prod_{j=1}^d
h_j$ and use
the notation $K_{h}(\cdot):=K ((\cdot-x_0)/h )/\Pi_h$ at
some points. Moreover, we define the neighborhood of $ x_0 $ of size $ h
$ as $ V_h:= \{x\in\bR^d \dvtx  K_{h}(x)>0 \} $ and assume for
simplicity that the kernel is chosen such that $ V_h\subseteq[0,1]^d$.
Next, we specify what we mean by a {contrast}:
%
\begin{definition}\label{def contrast}
A function $ \rho\dvtx  \bR\rightarrow[0,\infty)$ is called \textit{contrast}
(function) if it has the
following
properties:
\begin{enumerate}[3.]
\item$ \rho(\cdot) $ is convex, symmetric and $
\rho(0)=0 $;
\item the derivative $ \rho'(\cdot)$ of $ \rho(\cdot) $ is $ 1
$-Lipschitz and bounded;
\item the second derivative $ \rho''(\cdot) $ of $ \rho(\cdot) $ is
defined Lebesgue almost
everywhere and is $1$-Lipschitz with respect to the measure $\bP$.
Moreover, $ \|\rho''\|_\infty\leq1 $.
\end{enumerate}
\end{definition}

The constants in the Lipschitz condition and
the boundedness condition in the last definition are set
to $1$ for ease of exposition only. Well-known contrasts are the Huber
contrast (see \cite{Huber64}), for any scale $\gamma>0$ and $z\in\bR$,
%
\begin{equation}
\label{def huber contrast} \rho_{\mathrm{H},\gamma}(z):=\lleft\{ %
\begin{array} {l@{\qquad}l}
z^2/2, &\mbox{if } |z|\leq\gamma,
\\
\gamma\bigl(|z|-\gamma/2\bigr), &\mbox{otherwise}, \end{array} %
\rright.
\end{equation}
and the contrast induced by the arctan function (see \cite{Tsybakov82a})
%
\begin{equation}
\label{def arctan} \rho_{\op{arc},\gamma}(z):=\gamma z\op{arctan}(z/\gamma)-
\frac{\gamma^2}{2} \ln\bigl(1+z^2/\gamma^2\bigr).
\end{equation}
Note that the square loss and the absolute loss do not satisfy the above
definition. However, they can be mimicked by the Huber contrast with $
\gamma$ small (median) and $ \gamma$ large (mean). Let us define, for any
function $ \zeta$, the empirical measure as $ P_n\zeta:=\frac
{1}{n}\sum_{i=1}^n
\zeta(X_i,Y_i) $. We can now combine a kernel and a contrast to obtain
the {$ \lambda$-LPA estimator} $
\firstestimator\lambda({x_0}) $ of $ f^*(x_0) $ defined as:
%
\begin{eqnarray}
\label{def estimator with fixed lambda} &\firstestimator\lambda:=\arg\displaystyle \min_{f\in\cF}P_n
\lambda(f), \qquad& \mbox{where } \lambda(f) (x,y):= \rho \bigl(y-f(x) \bigr)
K_h(x)\nonumber \\[-8pt]\\[-8pt]
&&\mbox{ for }x\in[0,1]^d\mbox{ and }y\in\mr.\nonumber
\end{eqnarray}
The coefficients of the estimated polynomial can be
considered as estimators of the derivatives of the function
$ f^* $ at $x_0$. In this paper, however, we focus on the estimation of
$f^*(x_0)$.

%
%
%

\subsection{A first risk bound}

In this section, we present a risk bound for the estimator introduced
above. This estimator involves, in particular, fixed contrasts, kernels
and bandwidths.

To ease the presentation, we introduce some additional definitions.
First, we define
the best approximation of the target $
f^* $ in $ \cF
$ as
%
\begin{equation}
\label{coefficient_taylor} f^{0}:=\arg\min \Bigl\{\sup_{x\in
V_{h}}
\bigl|f(x)-f^*(x) \bigr| \dvtx  {f\in\cF},f(x_0)=f^*(x_0) \Bigr\}
\end{equation}
and the associated {bias term} as
%
\begin{equation}
\label{def bias term} b_h(\mathcal F):=\sup_{x\in V_h}\bigl|f^0(x)-f^*(x)\bigr|.
\end{equation}
The minimum is not necessarily unique, but all
minimizers work for our derivations. We then fix a multi-indexed vector
${t^0}=(t^0_{p_1,\dots,p_d})_{p\in\cP}$ such that $ \mathrm{P}_{t^0}=f^{0}
$. We recall that the entropy with bracketing of a set of functions
$\cA$ for a given radius $u>0$ with respect to a (pseudo)metric
$\Delta$
is the logarithm of the minimal number of pairs of functions
$(f_1^{(j)},f_2^{(j)})\in\cA\times\cA$ such that for any $f\in\cA$,
there is a couple $(f_1^{(j)},f_2^{(j)})$ such that $f_1^{(j)}\leq
f\leq
f_2^{(j)}$ and $\Delta(f_1^{(j)},f_2^{(j)})\leq u$. Here, in particular,
$ H_{\cF}(\cdot) $ denotes the entropy with bracketing of $ \cF$ with
respect to the\vspace*{1pt} pseudometric $\sqrt{\Pi_h\mathbb E
P_n [\lambda'(f_1)-\lambda'(f_2) ]^2}, f_1,f_2\in\cF$, where
%
\begin{equation}
\label{def first derivative lambda} \lambda'(f) (x,y):= \rho' \bigl(y-f(x)
\bigr) K_h(x) \qquad \mbox{for }x\in[0,1]^d\mbox{ and }y\in
\mr.
\end{equation}
The entropy $H_{\cF}(\cdot)$ cannot be calculated if the probability law
is unknown. However, it can be upper bounded invoking an upper bound
for the
pseudometric. For this, one may use that
\[
\sqrt{\Pi_h\mathbb EP_n \bigl[\lambda'(f_1)-
\lambda'(f_2) \bigr]^2}\leq \|K
\|_\infty\bigl\|t^{(1)}-t^{(2)}\bigr\|_1
\]
due to the continuity of $ \rho'(\cdot) $ and the definition of $ \cF$.
Here,
$ t^{(1)}$ and $ t^{(2)}\in[-M,M]^{|\cP|} $ are such that $ \mathrm{
P}_{t^{(1)}}=f_1 $ and $ \mathrm{P}_{t^{(2)}}=f_2 $, respectively. Therefore, the
entropy $H_{\cF}(\cdot)$ can be bounded by $ |\cP|$ times the
entropy of
$[-M,M]$ with respect to the Euclidean distance multiplied by
$\|K\|_\infty$ (this is, in particular, independent of $n$).

As a next step, we introduce the condition under which we derive the
risk bound.
%
\begin{condition}\label{condition consistency}
Let $\rho(\cdot)$ be a contrast, $K(\cdot)$ a kernel, $n\in\{
1,2,\ldots\}$,
and $h\in(0,1]^d$. We say that Condition \ref{condition consistency} is
satisfied if the smallest eigenvalue $\eigen$ of the matrix
\begin{eqnarray*}
\frac{1}{n}\sum_{i=1}^n\mathbb E
\biggl[U \biggl(\frac{X_i-{x_0}}{h} \biggr)U^\top \biggl(\frac
{X_i-{x_0}}{h}
\biggr)\rho'' \bigl(\sigma(X_i)
\xi_i \bigr) K_h(X) \biggr]
\end{eqnarray*}
is positive, ${n\Pi_h}\geq1$, and, defining
\begin{eqnarray*}
\delta_h:= \frac{2|\cP|^2}{\eigen} \biggl[\mathbb E\bigl[K_h(X)
\bigr]{b_h(\mathcal F)}+\frac{54\|\rho'\|_\infty (\sqrt{\me[\Pi_h
K^2_h(X)]}+\sfrac{\|K\|_\infty}{\sqrt{n\Pi_h}} )}{\sqrt{n\Pi
_h} ({{\ln(n|\cP|)}}+\int_0^{1}{H_{
\cF}^{1/2}(u)}\,\mathrm{d}u+{H_{\cF}(1)} )^{-1}} \biggr],
\end{eqnarray*}
that
%
\begin{equation}
\label{cond second derivative} 4\bigl(b_h(\mathcal F)+\delta_h\bigr)
\leq\inf_{x\in V_h}\frac{1}{n}\sum
_{i=1}^n\mathbb E \bigl[ \rho''
\bigl(\sigma(x)\xi_i \bigr) \bigr].
\end{equation}
\end{condition}

Condition \ref{condition consistency} can be interpreted in
the following sense: $n$ must be sufficiently large and $h$ appropriate
for the setting under consideration. In particular, $h$, as a function
of $n$, is usually chosen such that $h\to0$ and $n\Pi_h\to\infty$ as
$n\to\infty$ to satisfy Condition \ref{condition consistency}. We
postpone a detailed discussion of Condition \ref{condition consistency}
to after the main result of this section.

The variance term of the estimator is crucial for the following. To
state it explicitly, we need to introduce some more notation: First, we
introduce $ \lambda'' $ (similarly as $\lambda'$ in \eqref{def first
derivative lambda}) as
\[
\lambda''(f) (x,y):= \rho''
\bigl(y-f(x) \bigr) K_h(x) \qquad \mbox{for }x\in[0,1]^d\mbox{
and }y\in\mr.
\]
We then introduce the crucial quantity
%
\begin{equation}
\label{def Huber Variance} \mathrm{V}(\lambda):= \biggl(\frac{\sqrt{\Pi_h\mathbb E
P_n [\lambda'(f^*) ]^2}+
\|\rho'\|_\infty\|K\|_\infty\sfrac{\ln^2(n)}{\sqrt{n\Pi
_h}}}{\mathbb E
P_n\lambda''(f^*)} \biggr)^2.
\end{equation}
We call it nonasymptotic variance, since it plays the role of the
variance in the risk bounds in the theorems below. From
Condition \ref{condition consistency} and Definitions \ref{def kernel}
and \ref{def contrast}, we conclude that $ \mathrm{V}(\lambda)<\infty$. The
term $\|\rho'\|_\infty\|K\|_\infty\frac{\ln^2(n)}{\sqrt{n\Pi_h}}$
depends on $h$ and $n$. However, the bandwidth is typically chosen such
that $n\Pi_h\to\infty$ for $n\to\infty$ so that this term vanishes
asymptotically. Additionally, besides the normalization $\sqrt{\Pi_h}$
in front of the first term, a dependence on $h$ is given through
$\lambda$. We will discuss this after giving the main result of this
section. If $ h=(1,\ldots,1)^\top$ (parametric case), the
nonasymptotic variance $\mathrm{V}{(\lambda)}$ tends towards the asymptotic
variance $\op{AV}(\lambda)$ defined in \eqref{result normality huber}
as $n\to\infty$.

The main result of this section reads as the following.
%
\begin{theorem} \label{thr risk bound fix contrast}
Let $ \lambda$ be as in \eqref{def estimator with fixed lambda},
$n\in\{1,2,\dots\}$, and
$h\in(0,1]^d$ such that Condition \ref{condition consistency} is
satisfied. Then, for all $ q\geq1$,
\begin{eqnarray*}
&&\mathbb E \bigl|\hat{f}_{\lambda}({x_0} )-f^*(x_0)
\bigr|^q\\
&&\quad \leq C_q \biggl(b_h(\mathcal F)+
\biggl[27\int_0^{1}{H_{\cF}^{1/2}
(u )}\,\mathrm{d}u+\frac{4H_{\cF}(1)}{\ln^2(n)}+1 \biggr]\frac{\sqrt{\mathrm
{V}(\lambda)}}{\sqrt{n\Pi_h}} \biggr)^q+2^q
\frac{M^q}{n^2}
\end{eqnarray*}
for a constant $C_q$ ($C_q=4q|\cP|68^q\op{Gamma}(q)$ works, where $
\op{Gamma}(\cdot) $ is the classical Gamma function).
\end{theorem}

The proof can be easily deduced integrating the result of Proposition~\ref{prop uniform deviations} and using Proposition~\ref{prop large
deviation} (the propositions can be found in the \hyperref[app]{Appendix}, see also the
more detailed arXiv version).
%
\begin{remark}\label{remark th 1 constant}
In contrast to Huber's asymptotic results (see \cite{Huber64} and also
\cite{Arcones05,Tsybakov82a,Tsybakov82b,Tsybakov83,Tsybakov86}), the
above theorem holds for finite (but sufficiently large) sample sizes
$n$. We note that the desired variance term $\mathrm{V}(\lambda)$ is found
up to constants, which are of minor interest for this paper. Moreover,
a wide range of designs (including degenerate designs, e.g.) and
noise levels (including zero noise, e.g.) is covered. Let us
compare this result to \cite{Gaiffas05}: assume that $d=1$ and the
noise $ (\sigma(X_i)\xi_i)_i $ is identically and independently
normal distributed with variance $\sigma>0 $, and consider the local
Huber estimator with $ \rho(\cdot)=\rho_{\mathrm{H},\ln(n)}(\cdot) $
\eqref{def huber contrast}, where $ \gamma=\ln(n) $, and the
indicator kernel $ K(\cdot)=\1_{[ -1/2,1/2]}(\cdot) $. As we
mentioned above, the Huber estimator, with a large parameter $ \gamma
$, mimics the local least squares estimator. Indeed it holds
\[
\frac{\sqrt{\mathrm{V}(\lambda)}}{\sqrt{n\Pi_h}}\asymp\frac{\sigma
}{\sqrt{n\int_{x_0-h/2}^{x_0+h/2}\mu(x)\,\mathrm{d}x}}.
\]
The term on the right-hand side is the classical standard deviation of
the local least squares estimator. Theorem~\ref{thr risk bound fix
contrast} then implies the results of \cite{Gaiffas05}, Theorem~1 and Proposition~1, in the Gaussian case and extends them to heteroscedastic,
heavy-tailed noises.
\end{remark}
%
\begin{remark}\label{remark th 1}
While the above bound is -- to the best of our knowledge -- already a new
result, the final goal is to provide a specific $\lambda$ that
minimizes this bound since the second term $ 2^qM^q/n^2 $ is
neglectable and since the bias
term $b_h(\mathcal F)$ is independent of $\lambda$. However, the
bandwidth $h$, which accounts for the smoothness of the target
function, is included in ${\mathrm{V}(\lambda)}$. This makes simultaneous
D- and S-adaptation difficult. The specific dependences of the
numerator and the denominator on $h$ can be deduced from\vspace*{-1pt}
%
\begin{eqnarray}
\label{def numerator variance} \Pi_h\mathbb EP_n\bigl[
\lambda'\bigl(f^*\bigr)\bigr]^2&=&\Pi_h\int
\mu(x)K_h^2(x)\int \bigl[\rho' \bigl(
\sigma(x)z \bigr) \bigr]^2n^{-1}\sum
_i g_i (z )\,\mathrm{d}z \,\mathrm{d}x
\end{eqnarray}
and
\begin{eqnarray}
\label{def numerator variance 2} \mathbb EP_n\lambda''
\bigl(f^*\bigr)&=&\int\mu(x)K_h(x)\int \rho''
\bigl(\sigma(x)z \bigr)n^{-1}\sum_i
g_i (z )\,\mathrm{d}z \,\mathrm{d}x.
\end{eqnarray}
We study this in detail in the following section for three examples.
\end{remark}
\begin{discussion*}[of Condition \ref{condition consistency}]
The condition $\Phi_h>0$ is fulfilled in many examples. Indeed, with a
change of variables and by the definition of $ K_h(\cdot) $, we obtain
\begin{eqnarray*}
&&\frac{1}{n}\sum_{i=1}^n\mathbb E
\biggl[U \biggl(\frac{X_i-{x_0}}{h} \biggr)U^\top \biggl(\frac
{X_i-{x_0}}{h}
\biggr)\rho'' \bigl(\sigma(X_i)
\xi_i \bigr) K_h(X) \biggr]
\\
&&\quad  =\int U (x )U^\top (x )\mu(x_0+hx)K(x)\int
\rho'' \bigl(\sigma(x_0+hx)z \bigr)
\frac{1}{n}\sum_{i=1}^ng_i(z)\,\mathrm{d}z\,\mathrm{d}x.
\end{eqnarray*}
According to \cite{Tsybakov08}, Lemma~1.6, a sufficient condition for
$\Phi_h>0$ is thus that
%
\begin{equation}
\label{condition tsetse} \mu(x_0+hx)K(x)\int \rho''
\bigl(\sigma(x_0+hx)z \bigr) \frac{1}{n}\sum
_{i=1}^ng_i(z)\,\mathrm{d}z>0
\end{equation}
for all $x$ in some set in the kernel support with positive Lebesgue measure.
Recall that $\mu(x_0+h \cdot)$ is positive almost everywhere in the
support of $ K(\cdot) $ since $\mu(\cdot)$ vanishes only at finitely
many points. The condition $\Phi_h>0$ is thus fulfilled if
%
\begin{equation}
\label{eq asymp} \inf_{x\in V_h}\int \rho''
\bigl(\sigma(x)z \bigr) \frac{1}{n}\sum_{i=1}^ng_i(z)\,\mathrm{d}z>0.
\end{equation}
This condition is
satisfied, for example, for all densities $g_i(\cdot)$ and bounded
$\sigma(\cdot)$ if the contrast function is strictly convex. This holds
true for $\rho_{\op{arc},\gamma}(\cdot)$ (see \eqref{def
arctan}).
The Huber contrast $\rho_{\mathrm{H},\gamma}(\cdot)$ (see \eqref{def
huber contrast}), however, is strictly convex on the interval $( -
\gamma, \gamma)$ only. It holds that
$\rho_{\mathrm{H},\gamma}''(\cdot)=\1_{[ -\gamma, \gamma]}(\cdot)$;
therefore, the densities $g_i(\cdot)$ have to satisfy the additional
constraint
\[
\inf_{x\in V_h}\int\1_{[ -\gamma,
\gamma]} \bigl(\sigma(x)z \bigr)
\frac{1}{n}\sum_{i=1}^ng_i(z)\,\mathrm{d}z>0
\]
to ensure $\Phi_h>0$ in this case. If we assume,
for simplicity, that the noise level is constant $
\sigma(\cdot)\equiv\sigma>0 $, the last constraint simplifies to
$\int_{-\gamma/\sigma}^{\gamma/\sigma}\frac{1}{n}\sum_{i=1}^ng_i(z)\,\mathrm{d}z>0$. So even for
the Huber contrast with a fixed $ \gamma>0 $, the assumption $\Phi
_h>0$ is
weaker than the standard assumption in the literature of $g_i(\cdot)$
being positive and continuous in the origin for all $i\in\{1,2,\dots
,n\}$.

For the other crucial part of the condition, we first note that for
$h\to0$, the quantity on the right-hand side of \eqref{cond second
derivative} tends to a positive constant if $\sigma(\cdot)$ is
continuous in $x_0$. Indeed, since $\rho''$ is $\bP$-continuous, it
holds that
%
\begin{equation}
\inf_{x\in V_h}\frac{1}{n}\sum_{i=1}^n
\mathbb E \bigl[ \rho'' \bigl(\sigma(x)
\xi_i \bigr) \bigr]\longrightarrow\frac{1}{n}\sum
_{i=1}^n\mathbb E \bigl[ \rho''
\bigl(\sigma(x_0)\xi_i \bigr) \bigr]
\end{equation}
as $h\to0$. Similarly as above, the quantity on the right-hand side can
be lower bounded by a positive constant for many contrasts and noise
densities. We now give the rate for the quantity $\delta_h$. It holds that
\begin{eqnarray*}
\delta_h\asymp {\eigen^{-1}} \biggl[\mathbb E
\bigl[K_h(X)\bigr]{b_h(\mathcal F)}+\frac{\sqrt {\mathbb E\Pi_h
K_h^2(X)}}{\sqrt{n\Pi_h}}{{
\ln(n)}}+\frac{\|K\|_\infty}{{n\Pi
_h}}\ln(n) \biggr],
\end{eqnarray*}
which should be (cf.\ \eqref{cond second derivative}) bounded by a
constant. Here, ``$\asymp$''
indicates the asymptotic dependence on $n$. The above display corresponds
to a so-called \textit{bias--variance decomposition} up to a factor $ \ln(n)
$. We also
note that if $ f^* $ is continuous as assumed in the standard
literature, the bias term tends to zero as $ h\to0 $. In the literature,
one typically chooses some couple of positive constants $
(\alpha_1,\alpha_2) $, and $ \sigma(\cdot) $ and
some bandwidth $h=(h_1,\dots,h_d)$ such that
%
\begin{equation}
\label{example choice h} {n}^{-\alpha_1/d} \bigl({\ln(n)} \bigr)^{\alpha_2}\leq
h_j\leq \bigl({\ln(n)} \bigr)^{-1}\qquad  \mbox{for all }j=1,
\dots,d,
\end{equation}
where we assume that $n$ is sufficiently large such that the above
inequalities can hold.
For appropriate $ (\alpha_1,\alpha_2) $, Condition \ref{condition
consistency} is then satisfied
for $n$ sufficiently large in many examples.
\end{discussion*}
\begin{Example}\label{ex1}
If, for example, the design is uniform ($\mu(\cdot)\equiv1$) and the noise
level homoscedastic ($\sigma(\cdot)\equiv\sigma>0 $), it holds that
$\Phi_h\asymp \mathit{const}$ and thus $\delta_h\asymp b_h(\mathcal F) +
{\ln(n)}/{\sqrt{n\Pi_h}}$. Choosing a bandwidth $h=h_n$ as in
\eqref{example choice h} with $ \alpha_1=1 $ and $ \alpha_2=4 $,
Condition \eqref{cond second derivative} is satisfied for $n$
sufficiently large.
\end{Example}
\begin{Example}\label{ex2}
For degenerated designs, however, it is possible that $\Phi_h\to0$ as
$h\to
0$. For example, let $d=1$, the noise level be homoscedastic
($\sigma(\cdot)\equiv\sigma>0 $), and
%
\begin{equation}
\label{design number 2} \mu(\cdot)=\frac{s+1}{x_0^{s+1}+(1-x_0)^{s+1}}|\cdot-\,x_0|^s
\1 _{[0,1]}(\cdot)
\end{equation}
with $s> -1$ and $x_0\in[0,1]$ (see \cite{Gaiffas05}). The density
explodes (for $s<0$) or vanishes (for $s>0$) at $x_0$, so that one will
either have a lot or very little observations in the vicinity of $x_0$.
This is reflected in $\delta_h$ (recall that $d=1$ and thus $h\in(0,1]$):
\begin{eqnarray*}
\delta_h\asymp b_h(\mathcal F) +\frac{{{\ln(n)}}}{\sqrt{nh^{s+1}}}.
\end{eqnarray*}
So, similarly as above, one may choose a bandwidth like in
\eqref{example choice h} with $ \alpha_1=1/(s+1) $ and $
\alpha_2=4/(s+1) $.
\end{Example}

We finally note that the concrete form of Condition \ref{condition
consistency} is due to the application of deviation inequalities for
bounded empirical processes. Similarly, we could relax the
boundedness condition on the empirical processes involved to Bernstein
conditions (see, e.g., \cite{vdGeer11b}). This allows to incorporate
unbounded contrasts such as the least
squares contrast and the factor
$\|\rho'\|_\infty$ in Condition~\ref{condition consistency} should
be replaced
by the factor $\sqrt{\me[\rho'(\sigma(X)\xi)]^2}$, where $\xi$
would be a
sub-Gaussian random variable.

\subsection{S-minimax results}

In this section, we deduce some corollaries adapted to simple examples
from the above results.

To start, we recall the notion of S-minimaxity. To this end, let
$\tilde{f}(x_0)$ be an estimator of $f^*(x_0)$ and $\cS$ a set of
functions. For
any $q>0$, we define the \textit{maximal risk} of
$\tilde{f}$ and the \textit{S-minimax risk} for $x_0$ and $\cS$ as
%
\begin{equation}
\label{def maximalrisk} R_{n,q} [\tilde{f},\cS ] :=\sup_
{f^*\in\cS}
\mathbb{E} \bigl|\tilde{f} ({x_0})-f^*(x_0) \bigr|^q
\quad \mbox{and}\quad  R_{n,q} [\cS ]:=\inf_{\bar{f}}R_{n,q}
[\bar{f},\cS ],
\end{equation}
respectively. The infimum on the right-hand side is taken over all
estimators. We can now define the \textit{S-minimax rates of convergence}
and the
\textit{(asymptotic) S-minimax estimators}:
%
\begin{definition}
\label{def_minimax} A sequence $\phi_n$ is an
S-minimax rate of convergence, and the estimator $\hat{f}$ is
an (asymptotic) S-minimax estimator with respect to the set $\cS$ if
\[
\label{def_lower_bounds}0<\liminf_{n\to\infty}\phi^{-q}_n
R_{n,q} [\cS ] \leq\limsup_{n\to\infty}
\phi^{-q}_n R_{n,q} [\hat{f}, \cS ] <\infty.
\]
\end{definition}

We can give some simple examples for one dimensional target functions,
that is, $d=1$. We call $\bH_1(\beta, L,M)$ H\"{o}lder space, with
parameters $\beta, L,M>0$, the set of $ \lfloor\beta\rfloor$-times
differentiable functions $ f\dvtx [0,1]\to\mr$ such that $ \|f^{(j)}\|
_\infty\leq M\mbox{ for all }
j\in\{0,1,\dots,\lfloor\beta\rfloor\} $ and satisfied the H\"
{o}lder continuity $ |f^{(\lfloor\beta\rfloor)}(x)-f^{(\lfloor\beta
\rfloor)}(y)|\leq
L|x-y|^{\beta-\lfloor\beta\rfloor} \mbox{ for all } x,y\in[0,1]^d $.

The following corollary
can now be easily deduced from Theorem~\ref{thr risk bound fix contrast}.
%
\begin{corollary}\label{coro minimax}
Consider the model in Example~\ref{ex1}, that is, uniform design
($\mu(\cdot)\equiv1$) and homoscedastic noise level
($\sigma(\cdot)\equiv\sigma>0$). Let $\beta, L$ and $M$ be positive
parameters. Moreover, let $\firstestimator{\lambda}$ be defined as in
\eqref{def estimator with fixed lambda} with $m=\lfloor\beta\rfloor$,
$h\asymp n^{ - 1/(2\beta+1)}$,
$\rho(\cdot)=\rho_{\op{arc},1}(\cdot)$ as in \eqref{def arctan}, and
$K(\cdot):=\1_{[ -1/2,1/2]}(\cdot)$. Then, it holds that
\begin{eqnarray*}
\Pi_h\mathbb EP_n\bigl[\lambda'\bigl(f^*
\bigr)\bigr]^2&=& \frac{1}{n}\sum_{i=1}^n
\me \bigl[ \rho_{\op{arc},1}' (\sigma \xi_i )
\bigr]^2,
\\
\mathbb EP_n\lambda''\bigl(f^*\bigr)&=&
\frac{1}{n}\sum_{i=1}^n\me
\rho_{\op{arc},1}'' (\sigma\xi_i )
\end{eqnarray*}
 and
\begin{eqnarray*}
 \limsup_{n\to\infty} n^{q\beta/(2\beta+1)}&R_{n,q}
\bigl(\firstestimator{\lambda},\bH _1(\beta, L,M) \bigr)<\infty.
\end{eqnarray*}
\end{corollary}

The rate $n^{ -\beta/(2\beta+1)}$ is a standard S-minimax
rate in the context Gaussian noise (see \cite{Tsybakov08}, Chapter~2).
Here, however, this rate is achieved for a large class of noise
distributions.

Similarly, one can deduce the next corollary.
%
\begin{corollary}\label{coro minimax 2}
Consider the model in Example~\ref{ex2}, that is, a degenerate design as in
\eqref{design number 2} with $s> -1 $ and a homoscedastic noise
level ($\sigma(\cdot)\equiv\sigma>0$). Let $\beta, L,$ and $M$ be
positive parameters. Moreover, let $\firstestimator{\lambda}$ be defined
as in \eqref{def estimator with fixed lambda} with
$m=\lfloor\beta\rfloor$, $h\asymp n^{ - 1/(2\beta+s+1)}$,
$\rho(\cdot)=\rho_{\op{arc},1}(\cdot)$ as in \eqref{def arctan}, and
$K(\cdot):=\1_{\{ -1/2,1/2\}}(\cdot)$. Then, it holds that
\begin{eqnarray*}
\Pi_h\mathbb EP_n\bigl[\lambda'\bigl(f^*
\bigr)\bigr]^2&=&\frac
{h^s}{x_0^{s+1}+(1-x_0)^{s+1}}\frac{1}{n}\sum
_{i=1}^n \me \bigl[\rho_{\op{arc},1}'
(\sigma\xi_i ) \bigr]^2,
\\
\mathbb EP_n\lambda''\bigl(f^*\bigr)&=&
\frac{h^s}{x_0^{s+1}+(1-x_0)^{s+1}}\frac
{1}{n}\sum_{i=1}^n
\me \rho_{\op{arc},1}'' (\sigma\xi_i
)
\end{eqnarray*}
 and
\begin{eqnarray*}
 \limsup_{n\to\infty} n^{q\beta/(2\beta+s+1)}&R_{n,q}
\bigl(\firstestimator{\lambda},\bH _1(\beta, L,M) \bigr)<\infty.
\end{eqnarray*}
\end{corollary}

Thus, the rate $n^{ -\beta/(2\beta+s+1)}$ is achieved. This
rate is S-minimax in the nonparametric regression with homoscedastic
Gaussian noise (see \cite{Gaiffas05}).
{Note that we have only considered examples with homoscedastic noises
here. For heteroscedastic noises, the dependence on $h$ can be very
involved for some contrast functions (cf.\ equations \eqref{def numerator
variance} and \eqref{def numerator variance 2}). But, as highlighted
by the next example, this is not always the case.}
%
\begin{corollary}\label{coro minimax 3}
Consider the model \eqref{model} with $d=1$, a degenerate design as in
\eqref{design number 2} with $s> -1 $, a heteroscedastic noise
level $\sigma(\cdot)\equiv|\cdot-\,x_0|^\alpha, 0\leq\alpha\leq
s/2$, and a noise $(\xi_i)_i$ with finite variance. Let $\beta, L$
and $M$ be
positive parameters. Moreover, let $\firstestimator{\lambda}$ be defined
as in \eqref{def estimator with fixed lambda} with
$m=\lfloor\beta\rfloor$, $h\asymp n^{ - 1/(2\beta+s-2\alpha+1)}$,
$\rho(\cdot)=\rho_{\mathrm{H},\ln(n)}(\cdot)$ as in \eqref{def huber
contrast}, and
$K(\cdot):=\1_{\{ -1/2,1/2\}}(\cdot)$. Then, it holds that
\begin{eqnarray*}
\Pi_h\mathbb EP_n\bigl[\lambda'\bigl(f^*
\bigr)\bigr]^2\asymp h^{s+2\alpha},\qquad  \mathbb EP_n
\lambda''\bigl(f^*\bigr)\asymp{h^s}
\end{eqnarray*}
and
\[
\limsup_{n\to\infty} n^{q\beta/(2\beta+s-2\alpha+1)}R_{n,q} \bigl(
\firstestimator{\lambda },\bH_1(\beta, L,M) \bigr)<\infty.
\]
\end{corollary}

This result illustrates the effect of small noise levels on the
rate and the possible compensations to degenerate (unfavorable)
designs. In particular, if $\alpha=s/2$, we get
the standard minimax rate $n^{-\beta/(2\beta+1)}$ as in Corollary~\ref{coro minimax}. We assume that $ \alpha$ is smaller than $ s/2
$, since otherwise the noise level is very small and the bandwidth
chosen is
thus as small as possible. We also recall that the noise level is
assumed to be bounded so that we only consider the case $ \alpha\geq0 $.

\section{A D-adaptive estimator for fixed bandwidths}\label{section
Dadaptive fix bandwidth}

In this section, we discuss the selection of the combined function
$\lambda$, that is, of the
kernel and
the contrast. For this, we introduce an oracle that minimizes the
bound in Theorem~\ref{thr risk bound fix contrast} above and then
provide an estimator that mimics this oracle. This estimator is then
D-adaptive, that is, adaptive with respect to the noise and the design
distributions.

To this end, we first introduce $\Lambda:=\Upsilon\times\cK$ as the
set of
possible combined functions $\lambda$ as in \eqref{def estimator with fixed
lambda} for a given set of contrasts $\Upsilon$, a given set of kernels
$\cK$, and a fixed bandwidth $h\in(0,1]^d$. For example, one may
consider a
subset of the set of Huber functions indexed by the scale $\gamma>0$
as set
of contrasts
$
\Upsilon:=\{\rho_{\mathrm{H},\gamma}(\cdot)\dvtx \gamma>0\}.
$
An example for the set of kernels is the set of indicator functions
with different supports as
\begin{eqnarray*}
&&\cK:=\bigl\{\1_{S(u)}(\cdot)\dvtx u\in[ - 1/2,1/2]^d\bigr\}
\\
&&\quad \mbox{for }S(u):=[ -1/2+u_1, 1/2+u_1]\times\cdots
\times[ -1/2+u_d, 1/2+u_d].
\end{eqnarray*}
This contains, in particular, the symmetric indicator kernel $\1
_{S(0)}(\cdot)$.
In this section, the bandwidth $h$ is fixed so that the bias term $
b_h(\mathcal F) $ in Theorem~\ref{thr risk bound fix contrast} is of
minor importance; we then introduce the {oracle} as the minimizer of
the variance \eqref{def Huber Variance}
%
\begin{equation}
\label{def oracle lambda} \lambda^*:=\arg\min_{\lambda\in\Lambda}\mathrm{V}(\lambda).
\end{equation}
To mimic the oracle $\lambda^*$,
we propose the estimator $\widehat\lambda$
%
\begin{eqnarray}
\label{def adaptive lambda} &&\widehat\lambda:=\arg\min_{\lambda\in\Lambda}\widehat{\mathrm
{V}}(\lambda),\nonumber\\[-8pt]\\[-8pt]
&&\quad   \mbox{ where } \widehat{\mathrm{V}}(\lambda):= \biggl(
\frac{\sqrt {\Pi_hP_n [\lambda'(\hat f_{\lambda}) ]^2}+
\|\rho'\|_\infty\|K\|_\infty\sfrac{\ln^2(n)}{\sqrt{n\Pi_h}}}{P_n
\lambda''(\hat f_{\lambda})} \biggr)^2.\nonumber
\end{eqnarray}
Note that we estimate the target function $ f^* $ by $ \hat f_{\lambda
}$ and
$ \mathbb EP_n [\lambda' (f^*
) ]^2 $ and $ \mathbb EP_n\lambda'' (f^* )$ by their
empirical versions $ P_n [\lambda' (\hat
f_{\lambda}
) ]^2 $ and $P_n\lambda'' (\hat f_{\lambda}
)$, respectively. The explicit expressions for the numerator and the
denominator can be obtained using
\begin{eqnarray*}
P_n \bigl[\lambda'(\hat f_\lambda)
\bigr]^2&=&\frac{1}{n}\sum_{i=1}^nK_h^2(X_i)
\bigl[\rho' \bigl(Y_i-\hat f_\lambda(X_i)
\bigr) \bigr]^2 \quad \mbox{and}\\
P_n \lambda''(
\hat f_\lambda)&=&\frac{1}{n}\sum_{i=1}^nK_h(X_i)
\rho'' \bigl(Y_i-\hat
f_\lambda(X_i) \bigr).
\end{eqnarray*}
We now show that the estimator $\hat{f}_{\widehat\lambda}$ that
results from \eqref{def
estimator with fixed lambda} and \eqref{def adaptive lambda} performs
-- up
to\vspace*{-1pt} constants -- as well as the oracle $\hat{f}_{\lambda^*}$. For this,
we define $ H_{\cF\times\Lambda}(\cdot) $ as the entropy with
bracketing of
$\cF\times\Lambda$ with respect to the
(pseudo)metric
%
\begin{eqnarray}
\label{def distance entropy} &\sqrt{\Pi_h\mathbb EP_n \bigl[
\kappa(f_1,\lambda_1)-\kappa (f_2,
\lambda_2) \bigr]^2} \vee\sqrt{\Pi_h
\mathbb EP_n \bigl[\lambda_1''(f_1)-
\lambda_2''(f_2)
\bigr]^2}
\end{eqnarray}
for any $f_1,f_2\in\cF, \lambda_1,\lambda_2\in\Lambda$, where
$
\kappa(f,\lambda):={\lambda'(f)}/ (\sqrt{\Pi_h\mathbb
EP_n[\lambda'(f^*)]^2}+\|\rho'\|_\infty\|K\|_\infty\times\allowbreak \sfrac{\ln
^2(n)}{\sqrt{n\Pi_h}} ).
$
We compute in the \hyperref[app]{Appendix} a bound for this entropy for the set of Huer
contrasts indexed by the scale.

Before giving the main result of this section, we give the necessary
assumptions.
%
\begin{condition}\label{condition full consistency}
Let $\Lambda=\Upsilon\times\cK$ be a set of functions as in \eqref
{def estimator with fixed lambda} where $ \Upsilon$ is a set of
contrasts as in Definition~\ref{def contrast} and $ \cK$ is a set of
kernels as in Definition~\ref{def kernel}, $n\in\{1,2,\dots\}$, and
$h\in(0,1]^d$. We say that Condition \ref{condition full consistency}
is satisfied if the smallest eigenvalue $\eigen$ (defined in Condition
\ref{condition consistency}) is positive, ${n\Pi_h}\geq\ln^4(n)$,
and, defining for any $ \lambda\in\Lambda$
\begin{eqnarray*}
\delta_h^*(\lambda):= \frac{2|\cP|^2}{\eigen} \biggl[\mathbb E
\bigl[K_h(X)\bigr]{b_h(\mathcal F)}+\frac{54\|\rho'\|_\infty (\sqrt{\me
[\Pi_h K^2_h(X)]}+\sfrac{\|K\|_\infty}{\sqrt{n\Pi_h}}
)}{\sqrt{n\Pi_h} ({{\ln(n|\cP|)}}+\int_0^{1}{H_{
\cF\times\Lambda}^{1/2}(u)}\,\mathrm{d}u+{H_{\cF\times\Lambda}(1)}
)^{-1}}
\biggr],
\end{eqnarray*}
it holds for all $ {\lambda\in\Lambda} $
%
\begin{equation}
\label{cond star second derivative} 4\bigl(b_h(\mathcal F)+\delta_h^*(
\lambda)\bigr)\leq\inf_{x\in V_h}\frac
{1}{n}\sum
_{i=1}^n\mathbb E \bigl[ \rho''
\bigl(\sigma(x)\xi_i \bigr) \bigr].
\end{equation}
\end{condition}
%
\begin{condition}\label{condition variance estimation}
Additionally, we say that Condition \ref{condition variance
estimation} is satisfied if, defining
\begin{eqnarray*}
s_h(\lambda)&:=& \bigl(1\vee2\|K\|_\infty \bigr) \bigl[\delta
_h^*(\lambda)+b_h(\mathcal F) \bigr]\\
&&{}+27
\frac{ (1 \vee\|K\|
_\infty\|\rho'\|_\infty^{2} )}{\sqrt{n\Pi_h}} \biggl({{\ln \bigl(n|\cP|\bigr)}}+\int_0^{1}{H_{
\cF\times\Lambda}^{1/2}(u)}\,\mathrm{d}u+{H_{\cF\times\Lambda}(1)}
\biggr),
\end{eqnarray*}
it holds for all $ {\lambda\in\Lambda} $
%
\begin{equation}
\label{def radius 2} s_h(\lambda)\leq\frac{1}{2\|K\|_\infty}\min \bigl\{\mathbb
EP_n\lambda''\bigl(f^*\bigr),
\Pi_h \mathbb EP_n \bigl[\lambda'\bigl(f^*
\bigr) \bigr]^2 \bigr\}.
\end{equation}
\end{condition}

We discuss the above conditions after the following result.
%
\begin{theorem}\label{thr oracle bound}
Let $ \Lambda$ be a set of functions as in \eqref{def estimator with
fixed lambda}, $n\in\{1,2,\dots\}$, and
$h\in(0,1]^d$ such that Conditions \ref{condition full consistency}
and \ref{condition variance estimation} are satisfied. Then, for all $
q\geq1$,
\begin{eqnarray*}
&&\mathbb E \bigl|\hat{f}_{\widehat\lambda}({x_0} )-f^*(x_0)
\bigr|^q\\
&&\quad \leq T_q \biggl(b_h(\mathcal F)+
\biggl[27\int_0^{1}{H^{1/2}_{\cF\times
\Lambda}(u)}\,\mathrm{d}u+
\frac{4H_{\cF\times\Lambda}(1)}{\ln^2(n)}+1 \biggr]\frac{\sqrt{\mathrm{V}
(\lambda^*)}
}{\sqrt{n\Pi_h}} \biggr)^q+
\frac{5(2M)^q}{n^2}
\end{eqnarray*}
for a constant $T_q$ ($T_q=2q|\cP|117^q\op{Gamma}(q)$ works, where $
\op{Gamma}(\cdot) $ is the classical Gamma function).
\end{theorem}
%
\begin{remark}\label{remark D-adaptation}
Apart from the given assumptions, the estimator $\hat{f}_{\widehat
\lambda}(x_0)$ does not premise knowledge about the noise level
$\sigma(\cdot)$ and the densities $g_i(\cdot)$ and
$\mu(\cdot)$ but achieves -- up to constants -- the optimal variance
term $\mathrm{V}(\lambda^*) $ for all such functions. The estimator is
thus called D-adaptive
optimal (with respect to the set $\Lambda$). For example, for the Huber
contrast \eqref{def huber
contrast} indexed by the scale $ \gamma$, $
\Upsilon:=\{\rho_{\mathrm{H},\gamma}(\cdot)\dvtx \gamma>0\}
$, the estimator is D-adaptive minimax (Huber minimax)
for the set of contamination models, see Section \textup{A.1} in the arXiv
version. Finally, we mention that appropriate choices of the bandwidth
$h$ in the above result lead to S-minimax results.
\end{remark}
\begin{discussion*}[of Conditions \ref{condition full consistency}
and \ref{condition variance estimation}] Condition \ref{condition full
consistency} limits the possible sets of combined functions $\Lambda$ and
thus, in particular, the sets of possible contrast functions $\Upsilon
$. It
demands that all possible combined functions $\lambda\in\Lambda$ fulfill
Condition \ref{condition consistency}, which then leads to consistent
estimators (see Proposition~\ref{prop full consistency}) and to sets
of contrast with finite entropy. Condition \ref{condition full
consistency} demands, in
particular, that the right-hand side of \eqref{cond star second derivative}
is positive and, since the right-hand side of \eqref{cond star second
derivative} is upper bounded by $1$, that
$\sup_{\rho\in\Upsilon}\|\rho'\|_\infty$ does not increase too
rapidly with
$n$.

In the following, we illustrate these restrictions with an example. We
consider a homoscedastic model ($\sigma(\cdot)\equiv\sigma>0$) and
$\Upsilon$ equal to a set of Huber contrasts $\rho_{\mathrm{H},\gamma}(\cdot
)$ as in \eqref{def huber contrast} with scale parameter $\gamma\in
[\gamma_-,\gamma^+]$, $\gamma^+\geq\gamma_->0$. It holds that
$\sup_{\rho\in\Upsilon}\|\rho_{\mathrm{H},\gamma}'\|_\infty=\gamma^+$.
This implies that $\gamma^+$ must not increase too rapidly with $n$.
Moreover, it must hold that
\begin{eqnarray*}
\frac{1}{n}\sum_{i=1}^n\mathbb E
\bigl[ \rho''_{\mathrm{H},\gamma} (\sigma
\xi_i ) \bigr]=\int_{ - \gamma/\sigma}^{\gamma/\sigma}
\frac
{1}{n}\sum_{i=1}^ng_i(z)\,\mathrm{d}z
\geq\int_{ - \gamma_-/\sigma}^{\gamma
_-/\sigma}\frac{1}{n}\sum
_{i=1}^ng_i(z)\,\mathrm{d}z>0.
\end{eqnarray*}
For noise densities that are positive and continuous in the origin,
this condition is verified for all $\gamma_->0$. For more involved
noise densities (vanished at the origin), however, $\gamma_-$ has to
be chosen sufficiently large.

Condition \ref{condition variance estimation} is similar to Condition
\ref{condition full consistency} since $ s_h(\lambda)\asymp\delta
^*_h(\lambda) $. However, the terms in the minimum on the
right-hand side of \eqref{def radius 2}, can be small for a certain
design and noise level. The second term, vanishes if $\sigma(\cdot
)\equiv0$ since
$ \rho'(0)=0 $. Moreover, if the design degenerates (as in \eqref
{design number 2}) with a large $s$, $\mathbb EP_n\lambda''(f^*)$ and
$\Pi_h
\mathbb EP_n [\lambda'(f^*) ]^2$ then tend to zero faster
than $ s_h(\lambda) $ as $ n\to\infty$ (cf.\ \eqref{def numerator
variance} and \eqref{def numerator variance 2}). This is due to the
estimation of $\mathbb EP_n\lambda''(f^*)$ and $\Pi_h
\mathbb EP_n [\lambda'(f^*) ]^2$: if
$\sigma(\cdot)\equiv0$ or if the design degenerates, the above terms
are small (cf.\ \eqref{def Huber Variance}), and thus, the estimation
error of them (which is related to $ s_h(\cdot) $) obstructs their behavior.
\end{discussion*}
\section{{A D-adaptive and S-adaptive estimator}}\label{section ds
adaptation}
In this section, we introduce an estimator of $f^*(x_0)$ that is
simultaneously S- and
D-adaptive. For this, we apply the data-driven
procedure introduced above to select the contrast and the kernel and a
modification of the
data-driven Lepski's method to select the bandwidth. In the first part, we
consider isotropic, locally polynomial target functions, in the second part
anisotropic, locally constant functions. To simplify the exposition, we
present asymptotic results only.

The LPA is designed for functions that can be locally approximated by
polynomials. This is, for example, the case for H\"older classes, which we
define (similarly as
in \cite{Bertin04}) as
%
\begin{definition}\label{def holder space} Let $
\vec{\beta}:=(\beta_1,\ldots,\beta_d)\in\,]0,+\infty[^d $ such that
$
\lfloor\beta_1\rfloor=\cdots=\lfloor\beta_d\rfloor=:\lfloor
\beta\rfloor
$, and
let $ L,M>0 $. The function $ s\dvtx [0,1]^d\rightarrow[ -M,M] $
belongs to the
anisotropic H\"{o}lder
Class
$\mathbb{H}_d(\vec{\beta},L,M)$ if for all $
x,x_0\in[0,1]^d $
\begin{eqnarray*}
\bigl|s(x)-\mathrm{P}(s) (x-x_0)\bigr|&\leq& L\sum
_{j=1}^d|x_j-x_{0,j}|^{\beta_j}
\quad \mbox{and}
\\
 \sum_{p\in\cS_{\lfloor\beta\rfloor}}\sup_{x\in[0,1]^d}
\biggl\llvert \frac{\partial^{|
p|}s(x)}{\partial
x_1^{p_1}\,\cdots\,\partial x_d^{p_d}}\biggr\rrvert &\leq& M ,
\end{eqnarray*}
where $ \mathrm{P}(s)(x-x_0) $ is the Taylor polynomial
of $s$
of order $\lfloor\beta\rfloor$ at $ x_0 $, and $x_j$ and
$x_{0,j}$ are the $j$th components of $x$ and $x_0$, respectively.
\end{definition}

The parameter $ \vec\beta$ is usually unknown; thus, it is
desirable to have an estimator that is
adaptive with respect to $\vec\beta$. This motivates the following
definition, where $
\Psi:= \{\psi_n(\vec\beta) \}_{\vec\beta\in\cM} $ is
a given
family of normalizations for a set of parameters $ \cM$:
%
\begin{definition}\label{def S-adaptive}
The family $\Psi$ is called admissible if there exists an estimator
$\hat{f}_n$ such that
\begin{eqnarray*}
\limsup_{n\rightarrow\infty}\sup_{\vec\beta\in
\cM}
\psi^{-q}_n(\vec\beta) R_{n,q} \bigl(
\hat{f}_n,\mathbb{H}_d(\vec {\beta},L,M) \bigr)<\infty.
\end{eqnarray*}
The estimator $\hat{f}_n$ is then called $\Psi$-{adaptive} in the
S-minimax sense.
\end{definition}

We distinguish two cases in the following: First, we consider the special
case of isotropic H\"older classes, that is,
$\beta_1=\cdots=\beta_d$. These classes only require a common bandwidth
for all dimensions that is chosen with the standard version of Lepski's
method (see
\cite{Lepski90} and \cite{Lepski_Mammen_Spokoiny97}). Afterwards, we
allow for anisotropic H\"{o}lder classes. These classes necessitate a
separate bandwidth for
every dimension of the domain under consideration. The standard version
of Lepski's
method is not applicable in this case, because it requires a monotonous
bias. We
circumvent this problem using a
modified version of Lepski's method as described in
\cite{Kerkyacharian_Lepski_Picard01} and \cite{Lepski_Levit99}.

\subsection{A fully adaptive estimator for isotropic, locally polynomial
functions}\label{section adaptation iso}
Here, we consider isotropic H\"{o}lder classes with $ \beta\in(0,m+1]
$, where
$m$ is the degree of the estimator $ \hat f_\lambda$ and may be chosen
arbitrarily large. Therefore, only one bandwidth
$h_{\op{iso}}=h_1=\cdots=h_d>0$ has to be selected. Geometrically, this
means that we select a hypercube in $\mr^d$ with edge length $h_{\op{iso}}$
as domain of interest (in
contrast to the anisotropic case, where we select a hyperrectangle with
edge lengths $h_1,\dots,h_d$).

A major issue is the choice of the bandwidth. In the following, we assume
that the variance term
${\mathrm{V}(\lambda_{h_{\op{iso}}})}/({n h_{\op{iso}}^d})$ for
(see Definitions \eqref{def estimator with fixed lambda} and \eqref{def
Huber Variance})
\begin{eqnarray*}
\lambda_{h_{\op{iso}}}(f) (x,y):=\rho \bigl(y-f(x) \bigr) K_{h_{\op
{iso}}}(x)
\qquad \mbox{for all } x\in[0,1]^d,y\in\bR,
\end{eqnarray*}
is \textit{decreasing} in the bandwidth so that we can apply Lepski's
method. This
imposes an additional restriction on the design and the noise. After the
main result of this section, we give some examples for designs and noises
that fulfill this restriction. Next, we introduce the set of bandwidths
$ \cH^{\op{iso}}:=[h_{-},h^+] $, where
$0<h_-<h^+<1$ are defined as (cf.\ \eqref{example choice h})
%
\begin{equation}
\label{def bandwidth net} h_-:=\frac{\ln^{6/d}(n)}{n^{1/d}} \quad \mbox{and}\quad  h^+:=\frac{1}{\ln(n)}.
\end{equation}
Since the inequality $ h_-<h^+ $ has to be satisfied, $n$ is required
to be large enough.
We then introduce the isotropic M-estimator for any $
{h_{\op{iso}}}\in\cH^{\op{iso}} $ as
\begin{eqnarray*}
\label{def isotropic estimator} \hat f^{h_{\op{iso}}}_{\op{iso}}:=\arg\min
_{f\in\cF}P_n \widehat\lambda_{h_{\op{iso}}}(f),
\end{eqnarray*}
where
\[
\widehat\lambda_{h_{\op{iso}}}=\arg\min_{\lambda_{h_{\op{iso}}}\in
\Lambda} \widehat{\mathrm{V}}(
\lambda_{{h_{\op{iso}}}})
\]
and $ \widehat{\mathrm{V}}(\cdot) $ is defined in \eqref{def Huber Variance}.
Eventually, we introduce a net $
\netiso:=\{{h_{\op{iso}}}\in\cH^{\op{iso}}, \exists
m\in\bN\dvtx  {h_{\op{iso}}}=h^{+}\epsilon^m\}, \epsilon\in(0,1)
$, such that $ 1\leq|\netiso|\leq n $ and then
apply Lepski's method for isotropic functions (see \cite{Lepski90} and
\cite{Lepski_Mammen_Spokoiny97}) to
define the data-driven bandwidth $ \hat{h}_{\op{iso}} $:
%
\begin{eqnarray}
\label{def adaptive index} \hat{h}_{\op{iso}}&:=&\max \biggl\{{h_{\op{iso}}}\in\netiso\dvtx
\bigl\llvert \hat f^{{h_{\op{iso}}}}_{\op{iso}}(x_0)-\hat
f^{{h}_{\op{iso}}'}_{\op{iso}}(x_0)\bigr\rrvert \leq 15\sqrt2\bigl(B+
\penni\bigr)\sqrt\frac{{\widehat{\mathrm{V}}(\widehat\lambda_{{h_{\op{iso}}'}})}}{{n({h}_{\op{iso}}')^d}},
\nonumber
\\[-8pt]\\[-8pt]
&&\hphantom{\max \biggl\{}\mbox{ for all } {h}_{\op{iso}}'\in\netiso\mbox{ such that }
{h}_{\op{iso}}'\leq {h_{\op{iso}}} \biggr\},\nonumber
\end{eqnarray}
where $ \penni:=11\sqrt{\ln(n|\cH_\epsilon^{\op{iso}}|)} $ and $
B:=27\int_0^{1}{H^{1/2}_{\cF\times\Lambda}(u)}\,\mathrm{d}u+\frac{4H_{\cF
\times\Lambda}(1)}{\ln^2(n)}
$.

We now obtain on isotropic H\"{o}lder classes
$
\mathbb{H}^{\op{iso}}_d(\beta,L,M):=\mathbb{H}_d((\beta,\dots
,\beta),L,M)$,
 for all $\beta,L,M>0
$
the following result:
%
\begin{theorem}\label{Th adaptive iso}
Let $ \Lambda$ be a set of combined functions as in \eqref{def
estimator with fixed lambda} and $n\in\{1,2,\dots\}$ such that
Conditions \ref{condition full consistency} and \ref{condition variance
estimation} are satisfied for all $ h_{\op{iso}}\in\cH^{\op{iso}} $.
Then, for any $ x_0\in(0,1)^d$, any $\beta\in(0,m+1]$, and any $ L>0 $,
there exists a universal positive constant $ C>0 $ such that
\begin{eqnarray*}
{R_{n,q} \bigl[\hat f_{\op{iso}}^{\hat{h}_{\op{iso}}}(x_0),
\mathbb{H}^{\op
{iso}}_d({\beta},L,M) \bigr]} \leq C {\inf
_{{h_{\op{iso}}}\in\cH^{\op{iso}}} \biggl\{Ld{h}_{\op{iso}}^\beta+\penni\sqrt
\frac{{{{\mathrm{V}}(\lambda^*_{{h_{\op{iso}}}})}}}{{n{h}_{\op{iso}}^d}} \biggr\} ^{q}}\qquad \mbox{as }n\to\infty.
\end{eqnarray*}
\end{theorem}
%
\begin{remark}
This oracle inequality like result shows the simultaneous S- and
D-adaptation of the estimator. It generalizes results in
\cite{Cai_Zhou09}, which rely on the asymptotic equivalence of the block
median method, in two important aspects: First, it allows for
heteroscedastic regression models with random designs. Second, it does
not require that the noise densities are positive at their median and
thus allows for a wider range densities.
Finally, we note that
Lepski's method has been used for locally constant M-estimators in
\cite{Reiss_Rozenholc_Cuenod09} but -- to the best of our knowledge --
never to locally polynomial
M-estimators as it is done here.
\end{remark}
%
%
\begin{remark}
If only S-adaptation is considered, the conditions on $n$ can be
considerably\vspace*{-1pt} relaxed. In fact, assuming that $ \mathrm{V}(\lambda_{\cdot
}) $ is known, the estimator $ \hat f_{\lambda_{\tilde h_{\op{iso}}}}
$ of \eqref{def estimator with fixed lambda} can be applied\vspace*{-1pt} instead of
$\hat
f_{\op{iso}}^{\hat{h}_{\op{iso}}}$, where $ \tilde h_{\op{iso}} $ is
selected from \eqref{def adaptive index} replacing $ \widehat{\mathrm{V}}(\widehat\lambda_{{h_{\op{iso}}'}}) $ by
 ${\mathrm{V}}(\lambda_{{h_{\op{iso}}'}})$. The Conditions \ref{condition full
consistency} and
\ref{condition variance estimation} can then be replaced by Condition
\ref{condition consistency}.
\end{remark}
%
\begin{remark}
The variance term is decreasing for settings with indicator kernels and
homoscedastic noise
levels (as one can check easily starting from \eqref{def Huber
Variance}); for settings with indicator kernels, Huber contrasts, $
\sigma(\cdot)=1+|\cdot-\,x_0|^{\alpha} $ for a $\alpha\in[0,{\leq}
1/2]$, and $d=1$; and for many other settings. On
the contrary, the variance term can be increasing, for example, if the
noise level is symmetric in $
x_0 $ and convex.
\end{remark}
%
\begin{corollary}\label{coro isotropic}
Consider the model in Example~\ref{ex1} in the previous section with $
\mu(\cdot)\equiv1 $ (uniform design) and $ \sigma(\cdot)\equiv1 $
(homoscdastic noise level). For any $\beta\in(0,m+1]$ and any $ L>0 $,
it holds that
\begin{eqnarray*}
\limsup_{n\rightarrow\infty} \biggl(\frac{n}{\ln(n)} \biggr)^{q\beta/(2\beta+d)}
R_{n,q} \bigl[\hat f^{\hat
{h}_{\op{iso}}}_{\op{iso}}(x_0),
\mathbb{H}_d^{\op{iso}}({\beta },L,M) \bigr]<\infty.
\end{eqnarray*}
\end{corollary}

This corollary can be deduced minimizing the term on the right-hand
side of the last
theorem with a standard bias/variance trade-off.
%
\begin{remark}\label{remark iso minimax result}
The rate $
(\ln(n)/n )^{\beta/(2\beta+1)} $ in the
above corollary is admissible (cf.\ Definition~\ref{def S-adaptive}) over
isotropic H\"{o}lder spaces and is asymptotically
optimal (see \cite{Brow_Low96} and \cite{Lepski90}) up to the
logarithm $ \ln(n)
$, which is the usual price for the adapativity (see Section~\ref{section discussion} for more details). Moreover, the approach used to
deduce the above corollary presumes uniform designs and homescedastic
noises; however, more elaborate approaches, perhaps similar to the ones
in \cite{Gaiffas07}, may lead to comparable results for degenerate designs.
\end{remark}

\subsection{A fully adaptive estimator for anisotropic, locally constant
functions}\label{section full adaptive anisotropic}
In this part, we allow for anisotropic H\"older classes and
bandwidths. In return, we restrict ourselves to locally
constant functions, that is, $ m=0$ (and
thus $|\cP|=1$) and $ \cF=[-M,M] $. Moreover, we restrict
ourselves to
uniform designs ($\mu(\cdot)\equiv1$) and homoscedastic ($
\sigma(\cdot)\equiv\sigma\geq0 $) and identically distributed
noise ($ g_i(\cdot)\equiv g(\cdot) $ for
all $ i=1,\dots,n $). For this setting, we
introduce an S- and D-adaptive estimator of $f^*(x_0)$. The main properties
of this estimator are given in Theorem~\ref{Th adaptive aniso}.\\

We introduce an estimator for each bandwidth in the set $ \cH
:=[h_-,h^+]^d $, where $h_-$ and
$h^+$ are defined in the previous section. For this, we define
the variance term as
%
\begin{equation}
\label{def homo variance} \mathrm{V}(\rho,K):= \biggl(\frac{\sqrt{\int [\rho'(\sigma z)
]^2 g(z)\,\mathrm{d}z}+
\|\rho'\|_\infty\|K\|_\infty\sfrac{\ln^2(n)}{\sqrt {nh_{-}^d}}}{\int\rho''(\sigma
z) g(z)\,\mathrm{d}z} \biggr)^2,
\end{equation}
and
the oracle for a set of contrasts $ \Upsilon$ and a set of kernels $
\cK$ as
%
\begin{equation}
\label{def new oracle aniso} \bigl(\rho^*,K^*\bigr):=\arg\min_{\rho\in\Upsilon, K\in\cK}\mathrm{V}(
\rho,K).
\end{equation}
Next, we introduce an estimator of the variance term as
%
\begin{equation}
\label{def estimator variance aniso} \widehat{\mathrm{V}}(\rho,K):= \biggl(\frac{\sqrt{\sklfrac{1}{n}\sum_{i=1}^n [\rho'(Y_i-\hat
f_{\lambda_{h^+}}(X_i)) ]^2}+
\|\rho'\|_\infty\|K\|_\infty\sfrac{\ln^2(n)}{\sqrt {nh_{-}^d}}}{\sklfrac{1}{n}\sum_{i=1}^n\rho''(Y_i-\hat
f_{\lambda_{h^+}}(X_i))}
\biggr)^2,
\end{equation}
where $ \hat f_{\lambda_{h^+}} $ is defined in \eqref{def estimator
with fixed lambda} with $
\lambda=\lambda_{h^{+}}(f)(x,y):=\rho (y-f(x) ) K_{h^{+}}(x)
$, and
an estimator of the oracle as
%
\begin{equation}
\label{def choice contrast limit} (\hat\rho,\hat K):=\arg\min_{\rho\in\Upsilon, K\in\cK}\widehat{
\mathrm{V}}(\rho,K).
\end{equation}
We stress that the variance term $V$, the oracle $(\rho^*,K^*)$, and
their estimators
$\widehat V$ and $(\hat\rho,\hat
K)$ are independent of the bandwidth. We can finally introduce the
desired estimator $\hat f^h$ for all $h\in\cH$:
%
\begin{equation}
\label{def estimator adaptive lambda limit} \hat f^h:=\arg\min_{f\in\cF}n^{-1}
\sum_{i}\hat\rho \bigl(Y_i-f(X_i)
\bigr)\hat K_h(X_i).
\end{equation}
The crucial step is now the choice of the bandwidth with a modified version
of Lepski's method (see
\cite{Kerkyacharian_Lepski_Picard01} and \cite{Klutchnikoff05}).
First, we define for
all $ a,b\in\bR$ the scalar $a\vee b:=\max(a,b) $ and for all $
h,h'\in\cH$ the vector $ h\vee h':=(h_1\vee
h'_1,\ldots,h_d\vee h'_d) $. We then consider the two families of
Locally Constant Approximation (LCA) estimators (provoked by \eqref
{def estimator adaptive lambda limit})
\[
\bigl\{\hat f^h \bigr\}_{h\in\cH} \quad \mbox{and}\quad  \bigl\{\hat
f^{h,h'}:=\hat f^{h\vee h'} \bigr\}_{h,h'\in\cH^2}.
\]
Note that $ \hat f^{h,h'}=\hat f^{h',h} $ (commutativity). Similarly as
above, we then introduce a net $
\cH_\epsilon:=\{(h_{-},\dots,h_{-})\}\cup \{h\in\cH\dvtx  \forall
j=1,\ldots,d \ \exists
m_j\in\bN\dvtx  h_j=h^{+}\epsilon^{m_j} \}, \epsilon\in(0,1), $
such that
$|\cH_\epsilon|\leq n$ and set
$\penn:=11\sqrt{\ln(n|\cH_\epsilon|)}$. We finally select the bandwidth
according to
%
\begin{eqnarray}
\label{def adaptive aniso index} \hat h&:=&\max_{\preceq} \biggl\{h\in
\cH_\epsilon\dvtx  \bigl\llvert \hat f^{h,h'}(x_0)-\hat
f^{h'}(x_0)\bigr\rrvert \leq 16\bigl(B+\penn\bigr)\sqrt
\frac{{\widehat{\mathrm{V}}(\hat\rho,\hat
K)}}{{n\Pi_{h'}}}
\nonumber
\\[-8pt]\\[-8pt]
&&\hphantom{\max_{\preceq} \biggl\{}\mbox{for all } h'\in\cH_\epsilon\mbox{ such that }
h'\preceq h \biggr\}.\nonumber
\end{eqnarray}
The maximum is taken with respect to the order $ \preceq$ which we define
as\vspace*{-1pt}
$
h\preceq h' \ \Leftrightarrow\  \prod_{j=1}^dh_j\leq\prod_{j=1}^dh'_j.
$
Note, in particular, that the right-hand side of \eqref{def adaptive
aniso index} is
decreasing with respect to this order.

The above choice of the bandwidth leads to the estimator $\hat f^{\hat h}$
with the following properties:
%
\begin{theorem}\label{Th adaptive aniso}
Let $ \Lambda$ be a set of combined functions as in
\eqref{def estimator with fixed lambda} and let $n\in\{1,2,\dots\}$ such
that Conditions \ref{condition full consistency} and \ref{condition
variance estimation} are satisfied for all $ h\in\cH$. Then, for any $
x_0\in(0,1)^d$, any $\vec\beta\in(0,1]^d$, and any $ L>0 $, there exists
a universal constant $C$ such that
\begin{eqnarray*}
R_{n,q} \bigl[\hat f^{\hat
h}(x_0),
\mathbb{H}_d(\vec{\beta},L,M) \bigr]\leq C\inf_{h\in\cH}
\Biggl\{L\sum_{j=1}^dh_j^{\beta_j}+
\penn\sqrt \frac{{
{\mathrm{V}}(\rho^*,K^*)}}{{n\Pi_{h}}} \Biggr\}^q.
\end{eqnarray*}
\end{theorem}
We can also derive the following corollary from
Theorem~\ref{Th adaptive aniso} via a bias/variance trade-off.
%
\begin{corollary}\label{coro aniso}
For any $\vec\beta\in(0,1]^d$ and any $ L>0 $, it holds that
\begin{eqnarray*}
\limsup_{n\rightarrow\infty} \biggl(\frac{n}{\ln(n)} \biggr)^{q\bar
\beta/(2\bar\beta+1)}R_{n,q}
\bigl[\hat f^{\hat
h}(x_0),\mathbb{H}_d(\vec{
\beta},L,M) \bigr]<\infty,
\end{eqnarray*}
where $ \bar\beta:=  (\sum_j1/\beta_j  )^{-1}$ is the
harmonic average.
\end{corollary}
%
\begin{remark}
In contrast to the previous part, only locally constant functions are considered
here, which is due to the bias term (cf.\ Lemma~\ref{lem control bias
anisotrope}). To the best of our knowledge, the presented choice of the
bandwidth is the first application of
the \textit{anisotropic} Lepski's principle (\cite{Lepski_Levit99},
see also
\cite{Kerkyacharian_Lepski_Picard01,Klutchnikoff05,Goldenshluger_Lepski08}) for the selection of an anisotropic bandwidth for
nonlinear M-estimators. We also note that, comparing the adaptive rate $
(\ln(n)/n )^{\bar\beta/(2\bar\beta+1)} $ with the
optimal rate
in the white noise model (see \cite{Klutchnikoff05}), for example, one
finds that this rate is nearly optimal. We finally refer to the remarks
after Theorem~\ref{Th
adaptive iso}.
\end{remark}

\section{Discussion}\label{section discussion}


Let us detail on the assumptions and restrictions and highlight some
open problems:
\begin{enumerate}[6.]
\item
Instead of assuming that the densities $g_i(\cdot)$ are symmetric (cf.\
\cite{Huber64,Stone75}), it is sufficient that the sum $\sum_ig_i(\cdot)$ is symmetric. We are, however, not aware of examples
where this generalization is relevant.

\item The variance of the median estimator is $ 1/(4g^2(0)) $ which
implies a strong sensitive to the noise density at $0$. Moreover, the
estimation of $ g(0) $ (see \cite{Cai_Zhou09}, e.g.) requires
many observations near $f(x_0)$ in practice. On the contrary, Huber
contrast with scale $\gamma$ (allowed by our approach), the
denominator of the variance term \eqref{def Huber Variance} depends on
the mass of the noise density on the interval $ [-\gamma,\gamma
] $ instead of the mass at $0$.



\item To estimate the variance term \eqref{def Huber Variance}, we
plug an estimate $ Y_i-\hat f_{\lambda} $ of the residuals. Condition~\ref{condition full consistency} ensures the consistency of all
estimators in $ \Lambda$. This is considerably restrictive on the
initial family $\Lambda$.
This problem can be circumvented using a pre-estimator (e.g.,
with the contrast \eqref{def arctan}) instead of $ \hat f_{\lambda} $
for the estimation of the variance.

\item Lepski's method is very sensitive to outliers (see \cite
{Reiss_Rozenholc_Cuenod09}). To complement it with the adaptive
robustness of the estimator via the minimization of the variance term
can thus be interesting for many applications.




\item The variance term and its empirical version do not depend on the bias
term (see Theorem~\ref{thr risk bound fix contrast},
Definition \eqref{def adaptive lambda} and Remark~\ref{remark th 1})
and, more generally, not
on the specific model. The procedure presented in this
paper may thus be interesting for other models, such as high
dimensional settings (cf.\ \cite{Lambert_Zwald11}), for example.


\item The quantity $ 15\sqrt2(B+\penni)
$ in the threshold term in \eqref{def adaptive index} contains the
factor $\ln(n)$ and known but large constants. For applications, it
should usually be chosen considerably smaller (see \cite
{Lepski_Mammen_Spokoiny97}) and can probably be tuned with the \textit
{propagation} method \cite{Spokoiny_Vial09}, for example.

\item As mentioned in the \hyperref[intro]{Introduction}, Lepski-type procedures are also
useful to get S-adaptive confident bands (see \cite{Gine_Nickl10} and
references therein). This requires deviation inequalities that can be
derived along the presented lines (see Proposition~\ref{prop uniform
deviations}) but also a lower bound for the bias term of the estimator
(cf.\ \cite{Gine_Nickl10}, Condition 3, Section~3.2 and Section~3.5 for
Discussion), which seems not to be available here, since
robust M-estimators -- and thus the bias term -- do not have explicit
expressions. For our purposes, we circumvent this issue by using the
bias term of the criterions's derivative as an estimator of the
expected criterion's derivative, see Lemmas \ref{lem criterion
inversible} and \ref{lem control bias}. However, this way, we only
obtain an upper bound. We therefore suggest to establish first
S-adaptive confidence bands for the criterion's derivative viewed as an
estimator and then, using the smoothness of the contrast, confidence
bands with respect to a pointwise semi-norm or sup-norm.
%
%
%
\end{enumerate}

\section{Proofs of the main results}\label{section proofs main res}

Let us introduce some additional notation to simplify the
exposition. For this, we introduce
%
\begin{equation}
\label{def event consistence} \cF_{\delta}:= \bigl\{f=\mathrm{P}_t\in\cF\dvtx
\bigl\|t-t^0 \bigr\|_{\ell_1}\leq \delta \bigr\}
\end{equation}
as a ball in $ \cF$ with radius $\delta>0$ centered at $ f^0 $.
Furthermore, we denote the column vector of partial derivatives of the
criterion $P_n\lambda(\cdot)$ (defined in \eqref{def estimator with
fixed lambda}) by
%
\begin{eqnarray}
\label{Derivate lambda criterion} \tilde D_\lambda(\mathrm{P}_t):= \biggl(-
\frac{\partial
}{\partial
t_{p}} P_n\lambda(\mathrm{P}_t)
\biggr)_{p\in\cP} \qquad \mbox{for all } t\in\bR^{|\cP|},
\end{eqnarray}
and the ``parametric'' expectation with respect to the distribution
$\mathbb E^0$ of
$ (X,f^{0}(X)+\sigma(X)\xi)$ by
%
\begin{equation}
\label{def deterministic criterion} \mathbb E^0 \bigl[\tilde D_\lambda(\cdot)
\bigr].
\end{equation}
Next, for all $t\in\bR^{|\cP|}$, we introduce the \textit
{Jacobian matrix} $ J_D $ of $
\mathbb E^0 [\tilde
D_\lambda ] $ as
%
\begin{eqnarray}
\label{Jacobian_matrix} \bigl(J_D(\mathrm{P}_t)
\bigr)_{
p,
q\in\cP}:= \biggl(\frac{\partial}{\partial t_{q}}\mathbb E^0 \bigl[
\tilde D_\lambda^p(\mathrm{P}_t) \bigr]
\biggr)_{ p,
q\in\cP}= \biggl(\frac{\partial}{\partial
t_{q}}\mathbb E^0 \biggl[-
\frac{\partial}{\partial
t_{p}} P_n\lambda(\mathrm{P}_t) \biggr]
\biggr)_{p,q\in\cP},
\end{eqnarray}
where $ \tilde D_\lambda^{p}(\cdot) $ is the $ p $th component of $
\tilde D_\lambda(\cdot) $. The Jacobian matrix exists according to
Definition~\ref{def contrast} and
Fubini's theorem. Furthermore, the sup-norm on $\bR^{|\cP|}$ is
denoted by
$ \|\cdot\|_{\ell_\infty} $, and the vector of coefficients of the
estimated polynomial $ \hat
f_\lambda$ is denoted by $
\hat t_\lambda$. Moreover, we set
%
\begin{equation}
\label{def normalization term} c_\lambda:=\mathbb EP_n\lambda''
\bigl(f^*\bigr)
\end{equation}
and $b_h:=b_h(\mathcal F).$ We finally define $ \lambda'_\infty:=\|
\rho'\|_\infty\|K\|_\infty$ and for any $ z\geq0 $
%
\begin{eqnarray}
\label{def upper bound Massart} B_z:=27\int_0^{1}{H^{1/2}_{\cF\times\Lambda}(u)}\,\mathrm{d}u+
\frac{4H_{\cF
\times\Lambda}(1)}{\ln^2(n)}+7\sqrt{2z}+\frac{2z}{\ln^2(n)}.
\end{eqnarray}

\subsection{Auxiliary results}

The following propositions are basic for the proofs of the main
results. The proofs of the propositions are given in the \hyperref[app]{Appendix}.
%
\begin{proposition}\label{prop full consistency}
Let $\Lambda=\Upsilon\times\cK$ be a set of functions as in \eqref
{def estimator with fixed lambda} where $ \Upsilon$ is a set of
contrasts as in Definition~\ref{def contrast} and $ \cK$ is a set of
kernels as in Definition~\ref{def kernel}. Let $n\in\{1,2,\dots\}$
and $h\in(0,1]^d$ be such that Condition \ref{condition full
consistency} is satisfied. Then,
$
\bP (\bigcap_{\lambda\in\Lambda} \{\hat f_\lambda\in
\cF_{\delta_h^*(\lambda)} \} )\geq1-n^{-2},
$
where $ \delta_h^*(\cdot) $ is defined in Condition \ref{condition
full consistency}.
\end{proposition}

The following proposition allows us to control the deviations of the
process $ \tilde
D_\lambda(\cdot) $.
%
\begin{proposition}\label{prop large deviation}
For any $ z\geq0 $, it holds that
\[
\bP \biggl(\sup_{\lambda\in\Lambda}\sup_{f\in\cF_{\delta
_h^*(\lambda)}}
\frac{
\|\tilde
D_\lambda(f)-\mathbb E [\tilde D_\lambda(f) ]\|_{
\ell_\infty}}{\sqrt{\Pi_h\mathbb EP_n [\lambda'(f^*) ]^2}+
\lambda'_\infty\sfrac{\ln^2(n)}{\sqrt{n\Pi_h}}}\geq \frac{B_z}{\sqrt{n\Pi_h}} \biggr)\leq2|\cP|\exp(-z),
\]
where $ B_z $ is defined in \eqref{def upper bound Massart}.
\end{proposition}
This proposition is directly deduced from Massart's Inequality (see the
arXiv version for details).

\begin{proposition} \label{prop uniform deviations}
Let $\Lambda$ be a set of functions as in \eqref{def estimator with
fixed lambda},
$n\in\{1,2,\dots\}$, and $h\in(0,1]^d$ be such that Condition \ref
{condition full consistency} is satisfied. Then, for any $ z\geq0 $,
it holds that
\[
\bP \biggl( \biggl\{\sup_{\lambda\in\Lambda} \biggl[ \bigl|\hat{f}_{
\lambda}({x_0}
)-f^*(x_0) \bigr|-2\frac{\sqrt{\mathrm{V}(\lambda)} B_z}{\sqrt{n\Pi
_h}} \biggr]\geq 3b_h \biggr
\}\cap\bigcap_{\lambda\in\Lambda} \{\hat f_\lambda \in
\cF_{\delta_h^*(\lambda)} \} \biggr)\leq 2|\cP|\exp(-z).
\]
\end{proposition}
We note that the constants 2 and 3 can be replaced by $\mathrm{o}(1)$.
%
\begin{proposition}\label{prop control huber variance}
Let $\Lambda=\Upsilon\times\cK$ be a set of functions as in \eqref
{def estimator with fixed lambda}
where $ \Upsilon$ is a set of contrasts as in Definition~\ref{def
contrast} and $ \cK$ is a set of kernels as in Definition~\ref{def
kernel}. Let $n\in\{1,2,\dots\}$ and $h\in(0,1]^d$ be such that
Condition \ref{condition full consistency} is satisfied.

Then,
$
\bP (\Delta )
\geq1-5/n^2,
$
where $ \Delta:=\bigcap_{\lambda\in\Lambda} \{\sqrt{\widehat
{\mathrm{V}}
(\lambda)}\in [\frac{\sqrt2
}{3}\sqrt{{\mathrm{V}}(\lambda)},
\sqrt{6}\sqrt{{\mathrm{V}}(\lambda)} ] \} $.
\end{proposition}
We note that the constants $\sqrt2/3$ and $\sqrt{6}$ can be replaced
by $\mathrm{o}(1)$.

\subsection{Proof of Theorem \texorpdfstring{\protect\ref{thr oracle bound}}{2}}
First, we set
$\Delta:=\bigcap_{\lambda\in\Lambda} \{\sqrt{\widehat{\mathrm{V}}
(\lambda)}\in [\frac{\sqrt2
}{3}\sqrt{{\mathrm{V}}(\lambda)},
\sqrt6\sqrt{{\mathrm{V}}(\lambda)} ] \}.
$
Then, we observe that, since $ \hat f_{\hat\lambda}\in\cF$, $ \sup_{f\in\cF}|f(x_0)|\leq M $, and $ |f^*(x_0)|\leq M $, the
risk can
be bounded by
\begin{eqnarray*}
\mathbb E \bigl|\hat{f}_{\hat\lambda}({x_0} )-f^*(x_0)
\bigr|^q&=&\mathbb E \bigl|\hat{f}_{\hat\lambda}({x_0}
)-f^*(x_0) \bigr|^q\1_{\Delta}+\mathbb E \bigl|
\hat{f}_{\hat\lambda}({x_0} )-f^*(x_0) \bigr|^q
\1_{\Delta^c}
\\
&\leq&\mathbb E \bigl|\hat{f}_{\hat\lambda}({x_0} )-f^*(x_0)
\bigr|^q\1_{\Delta}+(2M)^q\bP\bigl({
\Delta^c}\bigr).
\end{eqnarray*}
Using Proposition~\ref{prop control huber variance}, Lemma~\ref{lemma
decomposition},
the last inequality, and simple
computations, we obtain
%
\begin{eqnarray}
\label{eq thr 2} \mathbb E \bigl|\hat{f}_{\hat\lambda}({x_0}
)-f^*(x_0) \bigr|^q&\leq&\mathbb E \bigl|\hat{f}_{\hat\lambda}({x_0}
)-f^*(x_0) \bigr|^q\1_{\Delta}+5(2M)^q/n^2
\nonumber
\\
&\leq&{2^{q}}\mathbb E{ \biggl( \bigl|\hat{f}_{
\hat\lambda}({x_0}
)-f^*(x_0) \bigr|-3b_h-\frac{6\sqrt{3\mathrm{V}(\lambda^*)}
B_0}{\sqrt{n\Pi_h}}
\biggr)_+^q\1_{\Delta}}
\\
&&{} +{2^{q}} \biggl(3b_h+\frac{6\sqrt{3\mathrm
{ V }
(\lambda^*)}
B_0}{\sqrt{n\Pi_h}}
\biggr)^q+5(2M)^q/n^2.\nonumber
\end{eqnarray}
Let us now bound the first term on the right-hand side of the last inequality.
To do so, we note that on the event $ \Delta$
%
\begin{equation}
\label{eq prop control huber variance} \sqrt{\mathrm{V}\bigl(\lambda^*\bigr)}\geq\sqrt{\frac{\widehat{\mathrm
{V}}(\lambda^*)}{6}}
\geq\sqrt{\frac{\widehat{\mathrm{V}}(\hat\lambda
)}{6}}\geq \sqrt{\frac{{\mathrm{V}}(\hat\lambda)}{27}}.
\end{equation}
%
Using the last inequality and integrating the result of Proposition~\ref{prop uniform deviations} with $ \varepsilon=10\sqrt{z} +\frac
{2z}{\ln^2(n)}$, we get (for more details see the arXiv version)
\begin{eqnarray*}
\mathbb E { \biggl( \bigl|\hat{f}_{\hat\lambda}({x_0}
)-f^*(x_0) \bigr|-3b_h-\frac{6\sqrt{3\mathrm{V}(\lambda^*)}
B_0}{\sqrt{n\Pi_h}}
\biggr)_+^q\1_{\Delta}} 
\leq T_q
\biggl(b_h+\frac{\sqrt{\mathrm{V}
(\lambda^*)}
B_0}{\sqrt{n\Pi_h}} \biggr)^q.
\end{eqnarray*}
From \eqref{eq thr 2} and the last inequality, the theorem can be deduced.

\subsection{Proof of Theorem \texorpdfstring{\protect\ref{Th adaptive iso}}{3}}
For ease of exposition, we set $ B_0=B $ (cf.\ \eqref{def upper bound
Massart}), $k:={h}_{\op{iso}}$, and $\hat k:=\hat{h}_{\op{iso}}$.
Then, one may verify that the \textit{oracle bandwidth}
\[
k^*:=\arg\min_{k\in\cH^{\op{iso}}} \biggl\{Ldk^\beta+c\bigl({B_0+
\penni }\bigr)\sqrt\frac{{{{\mathrm{ V}}(\lambda^*_k)}}}{{nk^d}} \biggr\}
\]
is well defined, where $ c $ is a constant chosen such that both terms
are equal at the point $ k^* $. Next, from Propositions \ref{prop full
consistency} and \ref{prop control huber variance} with $ h=(k,\ldots
,k) $, it
follows that
%
\begin{equation}
\label{eq thr iso 1 control event} \bP \bigl(\exists k\in\netiso,\exists\lambda_k\in
\Lambda\dvtx  \estimiso^k\notin\cF_{\delta_k^*(\lambda_k)} \bigr) \leq\sum
_{k\in\netiso}n^{-2}\leq n^{-1}
\end{equation}
and
%
\begin{eqnarray}
\label{eq thr iso 2 control event} \sum_{k\in\netiso}\bP \bigl(
\Delta^c_k \bigr) \leq\sum_{k\in\netiso}
\frac{5}{n^2}\leq5 n^{-1},
\end{eqnarray}
where $ \Delta_k:=\Delta$ is defined in Proposition~\ref{prop
control huber variance}.
Thus, we may restrict our considerations to the event $\bigcap_{k\in
\netiso,\lambda_k\in\Lambda} \{\estimiso^k\in\cF_{\delta
_k^*(\lambda_k)} \}\cap\Delta_k$, since we are only interested
in the asymptotic behavior. We now introduce $ k^*_\epsilon\in
\netiso $ such that
$ k^*_\epsilon\leq k^*\leq\epsilon^{-1} k^*_\epsilon$.

\subsubsection*{Control of the risk on the event $\{k^*_\epsilon\leq\hat
k\}$}
With the triangular inequality and Lemma~\ref{lemma decomposition}, we obtain
%
\begin{eqnarray}
\label{eq iso decomposition partie A}
&&\bigl|\estimiso^{\hat k}({x_0})-f^*({x_0})
\bigr|^q\1_{k^*_\epsilon
\leq\hat k}\nonumber\\[-8pt]\\[-8pt]
&&\quad \leq2^{q-1} \bigl( \bigl|
\estimiso^{\hat
k}({x_0})-\estimiso^{k^*_\epsilon}({x_0})
\bigr|^q\1_{k^*_\epsilon
\leq\hat k}+ \bigl|\estimiso^{k^*_\epsilon}({x_0})-f^*({x_0})
\bigr|^q \bigr).\nonumber
\end{eqnarray}
The first term on the right-hand side of the last inequality is
controlled using the procedure (\ref{def adaptive index}) to obtain
\[
\mathbb E \bigl[ \bigl|\estimiso^{\hat k}({x_0})-\estimiso
^{k^*_\epsilon}({x_0}) \bigr|^q\1_{k^*_\epsilon\leq\hat k} \bigr]\leq
\mathbb E \biggl[15\sqrt2\frac{\sqrt{{\widehat{\mathrm{V}}(\widehat
\lambda_{k^*_\epsilon})}}
(B_0+\penni)}{\sqrt{n(k^*_\epsilon)^d}} \biggr]^q.
\]
On the event $ \bigcap_{k\in\netiso}\Delta_k $, we get similarly as
in \eqref{eq prop control huber variance}
\[
\mathbb E \bigl[ \bigl|\estimiso^{\hat k}({x_0})-\estimiso
^{k^*_\epsilon}({x_0}) \bigr|^q\1_{k^*_\epsilon\leq\hat k} \bigr]\leq
\biggl(45\sqrt6\frac{\sqrt{{{\mathrm{ V}}(\lambda^*_{k^*_\epsilon})}}
(B_0+\penni)}{\sqrt{n(k^*_\epsilon)^d}} \biggr)^q.
\]
Recall that, by the definitions of the H\"{o}lder classes (Definition~\ref{def holder space}), we can control the bias for any $
\beta\in(0,m+1]$ and any $ k>0 $ by
%
\begin{equation}
\label{def control bias iso} b_k\leq\sup_{x\in
V_k}\bigl|\mathrm{P}
\bigl(f^*\bigr) (x-x_0)-f^*(x)\bigr|\leq Ldk^\beta,
\end{equation}
where $ \mathrm{P}(f^*)(x-x_0) $ is the Taylor Polynomial of $ f^* $
at $x_0$. So we can finally deduce from Theorem~\ref{thr oracle bound}
with $ h=(k,\ldots,k) $ and $ b_h=b_k $ a bound for the second term
in (\ref{eq iso decomposition partie A}) for $ n $ sufficiently large:
\[
\mathbb E \bigl|\estimiso^{k^*_\epsilon}({x_0})-f^*({x_0})
\bigr|^q\leq \cC_1 \biggl(Ld\bigl(k^*_\epsilon
\bigr)^{\beta}+\sqrt\frac{{{{\mathrm{V}}(\lambda
^*_{k^*_\epsilon})}}}{{n(k^*_\epsilon)^d}} \biggr)^q,
\]
where $ \cC_1 $ is a universal constant.
Using \eqref{eq iso decomposition partie A} and the above
inequalities, we have a control of the risk on the event $ \{
k^*_\epsilon\leq\hat k\} $:
%
\begin{eqnarray}
\label{eq control iso partie A} \mathbb E \bigl[ \bigl|\estimiso^{\hat k}({x_0})-f^*({x_0})
\bigr|^q\1 _{k^*_\epsilon\leq\hat k} \bigr]\leq\cC_2 \biggl(Ld
\bigl(k^*_\epsilon \bigr)^{\beta}+\bigl(B_0+\penni\bigr)\sqrt
\frac{{{{\mathrm{ V}}(\lambda^*_{k^*_\epsilon
})}}}{{n(k^*_\epsilon)^d}} \biggr)^q,
\end{eqnarray}
where $ \cC_1 $ is also a universal constant.

\subsubsection*{Control of the risk on the event $\{k^*_\epsilon>\hat k\}
$}
In order to control the risk on the
complementary
event, we observe that
%
\begin{eqnarray}
\label{eq iso decomposition partie Ac} \mathbb E \bigl[ \bigl|\estimiso^{\hat k}({x_0})-f^*({x_0})
\bigr|^q\1 _{k^*_\epsilon>\hat k} \bigr]\leq(2M)^{q}\bP
\bigl(k^*_\epsilon>\hat k\bigr).
\end{eqnarray}
We now show that the probability $ \bP(k^*_\epsilon>\hat k) $ is small.
According to the procedure \eqref{def adaptive index}, we have
\begin{eqnarray*}
\bP\bigl(k^*_\epsilon>\hat k\bigr)&\leq&\bP \biggl(\exists{k}'
\in\cH, {k}'< {k}^*_\epsilon\dvtx  \bigl\llvert
\estimiso^{k^*_\epsilon}(x_0)-\estimiso^{k'}(x_0)
\bigr\rrvert >15\sqrt2\frac{\sqrt{\widehat{\mathrm{V}}(\widehat\lambda_{k'})}
(B_0+\penni)}{\sqrt{n(k')^d}} \biggr)
\nonumber
\\
&\leq&2\sum_{k'\in\netiso: k'\leq k^*_\epsilon}\bP \biggl(\bigl\llvert
\estimiso^{k'}(x_0)-f^*(x_0)\bigr\rrvert >
\frac{15}{\sqrt2}\frac{\sqrt {\widehat{\mathrm{ V}}(\widehat\lambda_{k'})}(B_0+\penni)}{\sqrt {n(k')^d}} \biggr).
\end{eqnarray*}
On the event $ \bigcap_{k\in\netiso}\Delta_k $, we get similarly as
in \eqref{eq prop control huber variance}
\[
\bP\bigl(k^*_\epsilon>\hat k\bigr)\leq2\sum_{k'\in\netiso: k'\leq
k^*_\epsilon}
\bP \biggl(\bigl\llvert \estimiso^{k'}(x_0)-f^*(x_0)
\bigr\rrvert >5\frac{\sqrt{{\mathrm{V}}(\widehat\lambda_{k'})}(B_0+\penni)}{\sqrt {n(k')^d}} \biggr).
\]
Consequently,
%
\begin{eqnarray}
\label{eq th iso control complementary} \bP\bigl(k^*_\epsilon>\hat k\bigr)\leq2\sum
_{k'\in\netiso: k'\leq
k^*_\epsilon}\bP \biggl(\bigl\llvert \estimiso^{k'}(x_0)-f^*(x_0)
\bigr\rrvert >5\frac{\sqrt{{\mathrm{ V}}(\widehat\lambda_{k'})}(B_0+\penni)}{\sqrt {n(k')^d}} \biggr).
\end{eqnarray}
By definition, the oracle bandwidth $ {k}^* $ is the one which gives the
best trade-off. Thus, that the variance is decreasing, we obtain for
all $k'\leq k^*_\epsilon\leq k^*$
\begin{eqnarray*}
Ld\bigl(k'\bigr)^{\beta}&\leq& Ld\bigl(k^*_\epsilon
\bigr)^{\beta}\leq Ld\bigl(k^*\bigr)^{\beta}= \frac{\sqrt{{\mathrm{ V}}(\lambda^*_{k^*})}(B_0+\penni)}{\sqrt {n(k^*)^d}}
\\
&\leq&\frac{\sqrt{{\mathrm{ V}}(\lambda^*_{k^*_\epsilon})}(B_0+\penni
)}{\sqrt{n(k^*_\epsilon)^d}} \leq\frac{\sqrt{{\mathrm{ V}}(\lambda^*_{k'})}(B_0+\penni)}{\sqrt{n(k')^d}} \leq\frac{\sqrt{{\mathrm{ V}}(\widehat\lambda_{k'})}(B_0+\penni)}{\sqrt {n(k')^d}}.
\end{eqnarray*}
From \eqref{def control bias iso}, \eqref{eq th iso control
complementary}, and the last inequality, we get
\begin{eqnarray*}
\bP\bigl(k^*_\epsilon>\hat k\bigr)&\leq&2\sum
_{k
'\in\netiso: k'\leq k^*_\epsilon}\bP \biggl(\bigl\llvert \estimiso ^{k'}(x_0)-f^*(x_0)
\bigr\rrvert >2\frac{\sqrt{{\mathrm{ V}}(\widehat\lambda
_{k'})}(B_0+\penni)}{\sqrt{n(k')^d}}+3b_{k'} \biggr)
\nonumber
\\
&\leq&2\sum_{k'\in\netiso: k'\leq k^*_\epsilon}\bP \biggl(\sup
_{\lambda_{k'}\in\Lambda} \biggl[\bigl\llvert \estimiso ^{k'}(x_0)-f^*(x_0)
\bigr\rrvert -2\frac{\sqrt{{\mathrm{ V}}(\lambda
_{k'})}(B_0+\penni)}{\sqrt{n(k')^d}} \biggr]>3b_{k'} \biggr).
\end{eqnarray*}
Since $ \penni/\ln^2(n)\leq1 $ for $n$ sufficiently large, using the
definition of $ \penni$, Proposition~\ref{prop uniform deviations}
with $ h=(k',\ldots,k') $, $\lambda=\lambda_{k'}$, and $ z $ such
that $ B_z=(B_0+\penni) $, we obtain
\[
\bP\bigl(k^*_\epsilon>\hat k\bigr)\leq4|\cP|\sum
_{k'\in\netiso: k'\leq
k^*_\epsilon}\exp \biggl(-\frac{(\penni)^2}{100+4\penni/\ln
^2(n)} \biggr)\leq4|
\cP|n^{-1}.
\]
Then, in view of the last inequality, \eqref{eq thr iso 1 control
event}, \eqref{eq thr iso 2 control event}, \eqref{eq control iso
partie A} and \eqref{eq iso decomposition partie Ac}, we conclude that
\[
\mathbb E \bigl|\hat{f}^{\hat h}({x_0})-f^*({x_0})
\bigr|^q\leq\cC _2 \biggl(Ld\bigl(k^*_\epsilon
\bigr)^{\beta}+\bigl(B_0+\penni\bigr)\sqrt\frac{{{{\mathrm{
V}}(\lambda^*_{k^*_\epsilon})}}}{{n(k^*_\epsilon)^d}}
\biggr)^q\qquad \mbox {as }n\to\infty.
\]
By definition of $ k^* $ and $ k^*_\epsilon$ in the beginning of the
proof, the claim is proved.
\qed

\subsection{Proof of Theorem \texorpdfstring{\protect\ref{Th adaptive aniso}}{4}}
We set $ B=B_0 $. One may then verify that the \textit{oracle bandwidth}
\[
h^*:=\arg\min_{{h}\in\cH} \Biggl\{L\sum
_{j=1}^d\beta_j^{
-1}(h_j)^{\beta_j}+2
\frac{\sqrt{
{\mathrm{V}}(\rho^*,K^*)}(B_0+\penn)}{d\sqrt{n\Pi_{h}}} \Biggr\}
\]
is well defined. Define now the element $ h^*_{\epsilon} $ of $ \cH
_\epsilon$ such that
for all $ j=1,\ldots,d, h^*_{\epsilon,j}\leq h^*_{j}\leq\epsilon
^{-1}h^*_{\epsilon,j} $.
We then note that the estimator $ \hat
f^h $ is a constant function and $ f^0\equiv f^*(x_0) $, since we only consider
locally constant functions ($|\cP|=1$). To stress the importance of the
bandwidth, we set for any $ h\in\cH$,
$
\tilde\cD_h(\cdot):=\tilde D_{\tilde
\lambda_h}(\cdot)=n^{-1}\sum_{i}\hat\rho'(Y_i-\cdot)\hat
K_h(X_i)
$
and
%
\begin{equation}
\label{def re partial expectation anisotrope} \cD_h(\cdot):=\mathbb E \bigl[\tilde
D_{\tilde
\lambda_h}(\cdot) \bigr]=\int\hat K_h(x)\int\hat
\rho' \bigl(\sigma z+f^*(x)-\cdot \bigr)g(z)\,\mathrm{d}z\,\mathrm{d}x.
\end{equation}
Here, $\tilde\lambda_h(f)(x,y):=\hat\rho(y-f(x))\hat K_h(x)$ and $
(\hat\rho,\hat K) $ and $ \tilde D_{\lambda}(\cdot) $ are defined
in \eqref{def choice contrast limit} and \eqref{Derivate lambda
criterion}, respectively. Next, for uniform designs and homoscedastic noise
levels, the quantity $ c_{\lambda_h} $
%
\begin{equation}
\label{def re denominator} c_{\lambda_h}=c_{\rho}:=\int\rho''
(\sigma z)g(z)\,\mathrm{d}z,
\end{equation}
simplifies for any $ \lambda_h $ and does not depend on $h$.
Moreover, according to Lemma~\ref{lem aniso criterion inversible}, we have
for any $h\in\cH$, any $ \lambda\in\Lambda$, and any two constant
functions $ f,\tilde f\in\cF_{\delta_h^*(\lambda)} $
%
\begin{equation}
\label{eq thr 4 passage criterion} |f-\tilde f|\leq \tfrac{4}{3}c_{\hat\rho}^{-1}
\bigl|\cD_h(f)-\cD_h(\tilde f)\bigr|.
\end{equation}
Furthermore, from Propositions \ref{prop full consistency} and \ref
{prop control huber variance}, it
follows that
%
\begin{equation}
\label{eq thr 4 control event} \bP \bigl(\exists h\in\cH_\epsilon,\exists
\lambda_h\in\Lambda\dvtx  \hat f^h\notin\cF_{\delta_h^*(\lambda_h)}
\bigr) \leq\sum_{h\in\cH_\epsilon}n^{-2}\leq
n^{-1}
\end{equation}
and
%
\begin{eqnarray}
\label{eq thr 5 control event} \sum_{h\in\cH_\epsilon}\bP \bigl(
\Delta^c_h \bigr) \leq\sum_{h\in\cH_\epsilon}
\frac{5}{n^2}\leq5 n^{-1},
\end{eqnarray}
where $ \Delta_h:=\Delta$ is defined in Proposition~\ref{prop
control huber variance}.
Thus, we may restrict our considerations to the event $\bigcap_{h\in
\cH_\epsilon,\lambda_h\in\Lambda} \{\hat f^h\in\cF_{\delta
_h^*(\lambda_h)} \}\cap\Delta_h$, since we are only interested
on the asymptotic\vspace*{-1pt} behavior.
Moreover, we work on the event $
\cA:=\{h^*_\epsilon\preceq\hat h\} $ and its complement $\cA^c$
separately. For this, we decompose
the risk into $ R_\cA (\hat
f^{{h}},f^* ):=\mathbb E [ |\hat
{f}^{h}({x_0})-f^*({x_0}) |^q\1\{\cA\} ] $
and $ R_{\cA^c} (\hat
f^{{h}},f^* ):=\mathbb E [ |\hat
{f}^{h}({x_0})-f^*({x_0}) |^q\1\{\cA^c\} ] $.\\

\subsubsection*{Control of the risk on the event $\cA$}
With the triangular inequality and Lemma~\ref{lemma decomposition}, we obtain
%
\begin{eqnarray}
\label{eq decomposition partie A} R_\cA \bigl(\hat f^{\hat h},f^* \bigr)
\leq3^{q-1} \bigl[R_\cA \bigl(\hat f^{h^*_\epsilon,\hat h},\hat
f^{\hat
h} \bigr)+R_\cA \bigl(\hat f^{\hat h,h^*_\epsilon},\hat
f^{h^*_\epsilon} \bigr)+R_\cA \bigl(\hat f^{h^*_\epsilon},f^* \bigr)
\bigr].
\end{eqnarray}
Let us now control the first term on the right-hand side of the last
inequality. First, we observe that
%
\begin{equation}
\label{eq bound first term} R_\cA \bigl(\hat f^{h^*_\epsilon,\hat h},\hat
f^{\hat h} \bigr)\leq \mathbb E\sup_{h\in\cH: h\succeq h^*_\epsilon} \bigl|\hat
f^{h^*_\epsilon,h}({x_0})-\hat f^{h}({x_0})
\bigr|^q.
\end{equation}
%
Using \eqref{eq thr 4 passage criterion} and taking $ f=\hat
f^{h^*_\epsilon,h} $ and $ \tilde f=\hat f^{h} $, we then have
\begin{eqnarray*}
\bigl|\hat f^{h^*_\epsilon,h}({x_0})-\hat f^{h}({x_0})
\bigr| \leq2 c_{\hat\rho}^{-1} \bigl\llvert \cD_{h} \bigl(
\hat f^{h^*_\epsilon,h} \bigr)-\cD_{h} \bigl(\hat f^{h} \bigr)
\bigr\rrvert .
\end{eqnarray*}
Recall that, by definition, $ \tilde\cD_h(\hat f^h)=0 $ for all $
h\in\cH$. We then
obtain from the last inequality for any $ h\in\cH$
%
\begin{eqnarray}
\label{eq game m-estimation} &&\bigl|\hat f^{h^*_\epsilon,h}({x_0})-\hat
f^{h}({x_0}) \bigr| \nonumber\\
&&\quad \leq2 c_{\hat\rho}^{
-1}
\bigl(\bigl\llvert \cD_{h} \bigl(\hat f^{h^*_\epsilon,h} \bigr)-\cD
_{h^*_\epsilon\vee h} \bigl(\hat f^{h^*_\epsilon,h} \bigr)\bigr\rrvert
\\
&&\qquad \hphantom{2 c_{\hat\rho}^{
-1}
\big(}{} +\bigl\llvert \cD_{h^*_\epsilon\vee h} \bigl(\hat f^{h^*_\epsilon,h} \bigr)-
\tilde\cD_{h^*_\epsilon\vee h} \bigl(\hat f^{h^*_\epsilon,h} \bigr)\bigr\rrvert +\bigl
\llvert \tilde\cD_{h} \bigl(\hat f^{h} \bigr)-
\cD_{h} \bigl(\hat f^{h} \bigr)\bigr\rrvert \bigr).\nonumber
\end{eqnarray}
Denote by $\hat\lambda_h(f)(x,y)=\hat\rho (y-f(x) ) \hat
K_h(x)$ and $\tilde\delta_h:=\delta_h^*(\hat\lambda_h)\vee\delta
_{h\vee h_\epsilon^*}^*(\hat\lambda_{h\vee h_\epsilon^*})$, using
the last inequality and \eqref{eq bound first term}, we have
\begin{eqnarray*}
R_\cA \bigl(\hat f^{\hat h,h^*_\epsilon},\hat f^{\hat{h}} \bigr)
&\leq&2^{q-1} \mathbb Ec_{\hat\rho}^{
-q}\sup
_{h\in\cH_\epsilon}\sup_{f\in\cF_{\tilde\delta
_h}}2^q\bigl\llvert
\cD_{h}(f)-\cD_{h^*_\epsilon\vee h} (f)\bigr\rrvert ^q
\\
&&{}+2^{q}2^q\mathbb Ec_{\hat\rho}^{-q}\sup
_{h\in\cH: h\succeq
h^*_\epsilon}\sup_{f\in\cF_{\tilde\delta_h}}\bigl\llvert \tilde
\cD_{h}(f)- \cD_{h}(f)\bigr\rrvert ^q.
\end{eqnarray*}
Using Lemma~\ref{lem deviation criterion anisotrope} and Lemma~\ref
{lem control bias anisotrope} with $ h'=h^*_\epsilon$, there exists a
universal positive constant $ \cC$ such that
%
\begin{equation}
\label{eq th adaptive first term} R_\cA \bigl(\hat f^{\hat h,h^*_\epsilon},\hat
f^{\hat{h}} \bigr)\leq\cC \Biggl(L\sum_{j=1}^d
\bigl(h_{\epsilon,j}^*\bigr)^{\beta_j}+\frac
{\sqrt{\mathrm{V}(\rho^*,K^*)}(B_0+\penn)}{\sqrt{n\Pi_{h^*_\epsilon
}}}
\Biggr)^q.
\end{equation}
The second term on the right-hand side of \eqref{eq decomposition
partie A} is controlled by the procedure
(\ref{def adaptive aniso index}), which implies
\[
R_\cA \bigl(\hat f^{\hat h,h^*_\epsilon},\hat f^{h^*_\epsilon} \bigr)\leq
\mathbb E \biggl[16\frac{\sqrt{\widehat{\mathrm{V}}(\hat\rho,\hat K)}
(B_0+\penn)}{\sqrt{n\Pi_{h^*_\epsilon}}} \biggr]^q\1_\cA.
\]
On the event $ \bigcap_{h\in\cH_\epsilon}\Delta_h $,
%
\begin{eqnarray}
\label{eq th adaptive second term} R_\cA \bigl(\hat f^{\hat h,h^*_\epsilon},\hat
f^{h^*_\epsilon} \bigr)\leq \biggl(16\sqrt{6}\displaystyle \frac{\sqrt{{\mathrm{V}}(\rho^*,K^*)}
(B_0+\penn)}{\sqrt{n\Pi_{h^*_\epsilon}}}
\biggr)^q.
\end{eqnarray}
By the definition of the H\"{o}lder class (Definition~\ref{def holder
space}) and $ b_h $ (Definition \eqref{def bias term}), we can\vspace*{1pt} control
the bias for any $ h\in\cH$:
$
b_h\leq\sup_{x\in
V_h}|f^*(x_0)-f^*(x)|\leq
L\sum_{j=1}^dh_j^{\beta_j}.
$
Finally, with Theorem~\ref{thr oracle bound}, we can bound the third term
in (\ref{eq decomposition partie A}): There exists a universal
positive constant $ \cC$ such that
\[
R_\cA \bigl(\hat f^{h^*_\epsilon},f^* \bigr)\asymp\cC \Biggl(L\sum
_{j=1}^d\bigl(h_{\epsilon,j}^*
\bigr)^{\beta_j}+\frac{\sqrt{{\mathrm{V}}(\rho^*,K^*)}
B_0}{\sqrt{n\Pi_{h^*_\epsilon}}} \Biggr)^q.
\]
Using \eqref{eq decomposition partie A}, \eqref{eq th adaptive first
term}, \eqref{eq th adaptive second term}, and the last inequality, we
have a control of the risk on the event $ \cA$ such that
%
\begin{eqnarray}
\label{eq control partie A} R_\cA \bigl(\hat f^{\hat h},f^* \bigr)\leq\cC
\Biggl(L\sum_{j=1}^d\bigl(h_{\epsilon,j}^*
\bigr)^{\beta_j}+\frac{\sqrt{\mathrm{V}(\rho
^*,K^*)}(B_0+\penn)}{\sqrt{n\Pi_{h^*_\epsilon}}} \Biggr)^q
\end{eqnarray}
as $ n\to\infty$ and for a universal positive constant $ \cC$.

\subsubsection*{Control of the risk on the event $\cA^c$} In order to
control the risk on the
complementary event $\cA^c$, we observe that
%
\begin{eqnarray}
\label{eq cauchy Sharwz partie Ac} R_{\cA^c} \bigl(\hat f^{\hat h},f^* \bigr)
\leq(2M)^{q}\bP\bigl(\cA^c\bigr).
\end{eqnarray}
We now show that the probability $
\bP(\cA^c)
$ is small.
According to the construction of the procedure \eqref{def adaptive
aniso index}, the event $\cA^c$ implies that
there exists a
$ {h}'\in\cH_\epsilon$ such that $ {h}'\preceq h^*_\epsilon$ and
\[
\bigl\llvert \hat f^{h^*_\epsilon,{h}'}(x_0)-\hat f^{{h}'}(x_0)
\bigr\rrvert >16\frac{\sqrt {\widehat{\mathrm{V}}(\hat\rho,\hat K)}(B_0+\penn)}{\sqrt{n\Pi_{h'}}}.
\]
Using \eqref{eq thr 4 passage criterion} and taking $ f=\hat
f^{h^*_\epsilon,h'} $ and $ \tilde f=f^{h'} $, we have on the event $
\cA^c $
\[
\frac{4} 3 c_{\hat\rho}^{-1} \bigl\llvert
\cD_{h'} \bigl(\hat f^{h^*_\epsilon,h'} \bigr)-\cD_{h'} \bigl(
\hat f^{h'} \bigr)\bigr\rrvert >16\frac{\sqrt{\widehat{\mathrm{V}}(\hat\rho,\hat K)}(B_0+\penn
)}{\sqrt{n\Pi_{h'}}}.
\]
From the last inequality, we obtain (cf.\ \eqref{eq game m-estimation})
\begin{eqnarray*}
&&\frac{4} 3 c_{\hat\rho}^{
-1}\sup_{f\in\cF_{\tilde\delta_{h'}}}
\bigl\llvert \cD_{h'}(f)-\cD _{h^*_\epsilon\vee h'} (f)\bigr\rrvert +
\frac{8} 3 c_{\hat\rho}^{-1}\sup_{f\in\cF_{\tilde
\delta_{h'}}}
\bigl\llvert \tilde\cD_{h'}(f)- \cD_{h'}(f)\bigr\rrvert\\
&&\quad  >16
\frac{\sqrt{\widehat{\mathrm{V}}(\hat\rho,\hat
K)}(B_0+\penn)}{\sqrt{n\Pi_{h'}}}.
\end{eqnarray*}
Together with Lemma~\ref{lem control bias anisotrope}, this yields
\[
\frac{5} 3L\sum_{j=1}^d
\bigl(h_{\epsilon,j}^*\bigr)^{\beta_j}+\frac{8} 3
c_{\hat
\rho}^{-1}\sup_{f\in\cF_{\tilde\delta_{h'}}}\bigl\llvert \tilde
\cD_{h'}(f)- \cD_{h'}(f)\bigr\rrvert >16\frac{\sqrt{\widehat{\mathrm{V}}(\hat\rho,\hat
K)}(B_0+\penn)}{\sqrt{n\Pi_{h'}}}.
\]
On the event $ \bigcap_{h\in\cH_\epsilon}\Delta_h $, we get
similarly as in \eqref{eq prop control huber variance}
\[
\frac{5} 3L\sum_{j=1}^d
\bigl(h_{\epsilon,j}^*\bigr)^{\beta_j}+\frac{8} 3
c_{\hat
\rho}^{-1}\sup_{f\in\cF_{\tilde\delta_{h'}}}\bigl\llvert \tilde
\cD_{h'}(f)- \cD_{h'}(f)\bigr\rrvert >\frac{16\sqrt2} 3
\frac{\sqrt{{\mathrm{V}}(\hat
\rho,\hat K)}(B_0+\penn)}{\sqrt{n\Pi_{h'}}},
\]
this implies
\[
c_{\hat\rho}^{-1}\sup_{f\in\cF_{\tilde\delta_{h'}}}\bigl\llvert \tilde
\cD_{h'}(f)-\cD_{h'}(f)\bigr\rrvert >\frac{16\sqrt2}{8}
\frac{\sqrt{{\mathrm{V}}(\hat\rho,\hat
K)}(B_0+\penn)}{\sqrt{n\Pi_{h'}}}-\frac{5} 8 L\sum_{j=1}^d
\bigl(h_{\epsilon,j}^*\bigr)^{\beta_j}.
\]
By definition, the oracle bandwidth $ h^*_\epsilon$ is the one which
gives the
best trade-off. Thus by definition of $ h^*_\epsilon$, for all
$h'\preceq h^*_\epsilon\preceq h^*$
\begin{eqnarray*}
L\sum_{j=1}^d\bigl(h_{\epsilon,j}^*
\bigr)^{\beta_j}&\leq& L\sum_{j=1}^d
\bigl(h_{j}^*\bigr)^{\beta_j}= \frac{\sqrt{{\mathrm{V}}(\rho^*,K^*)}(B_0+\penn)}{\sqrt{n\Pi
_{h^*}}}
\\
& \leq&\frac{\sqrt{{\mathrm{V}}(\rho^*,K^*)}(B_0+\penn)}{\sqrt{n\Pi
_{h^*_\epsilon}}} \leq\frac{\sqrt{{\mathrm{V}}(\rho^*,K^*)}(B_0+\penn)}{\sqrt{n\Pi_{h'}}} \\
&\leq&\frac{\sqrt{{\mathrm{V}}(\hat\rho,\hat K)}(B_0+\penn)}{\sqrt {n\Pi_{h'}}}.
\end{eqnarray*}
From the last two inequalities, we obtain on the event $ \cA^c $
\[
c_{\hat\rho}^{-1}\sup_{f\in\cF_{\tilde\delta_h}}\bigl\llvert \tilde
\cD_{h'}(f)-\cD_{h'}(f)\bigr\rrvert >\frac{\sqrt{{\mathrm{V}}(\hat\rho,\hat K)}(B_0+\penn)}{\sqrt{n\Pi_{h'}}}.
\]
Then, we have a control of the following probability
\[
\bP\bigl(\cA^c\bigr) \leq\sum_{h'\in\cH_\epsilon: h'\preceq h^*_\epsilon}
\bP \biggl(\sup_{\rho,K}\sup_{f\in\cF_{\tilde\delta_h}}
\frac
{\llvert \tilde\cD_{h'}(f)-\cD_{h'}(f)\rrvert }{c_{\rho}\sqrt{{\mathrm
{V}}(\rho,K)}}>\frac{B_0+\penn}{\sqrt{n\Pi_{h'}}} \biggr).
\]
Using $ \penn/\ln^2(n)\leq1 $ and Propostion \ref{prop large
deviation} with $z$ such that $ B_z=B_0+\penn$, we deduce that
\begin{eqnarray*}
\bP\bigl(\cA^c\bigr) &\leq\displaystyle \sum
_{h'\in\cH_\epsilon: h'\preceq h^*_\epsilon}\exp \biggl(-\displaystyle \frac{(\penn)^2}{100+4\penn/\ln^2(n)} \biggr)\leq
n^{-1}.
\end{eqnarray*}
From \eqref{eq cauchy Sharwz partie Ac} and the last inequality, we
obtain on the event $ \cA^c $:
$
R_{\cA^c} (\hat f^{\hat h},f^* )\leq(2M)^qn^{-1}.
$
Then, in view of the last inequality, \eqref{eq thr 4 control event},
\eqref{eq thr 5 control event} and \eqref{eq control partie A}, we
conclude that there exists a universal positive constant $ \cC$ such that
\begin{eqnarray*}
\mathbb E \bigl|\hat{f}^{\hat h}({x_0})-f^*({x_0})
\bigr|^q\leq\cC \Biggl(L\sum_{j=1}^d
\bigl(h_{\epsilon,j}^*\bigr)^{\beta_j}+\frac{\sqrt{\mathrm
{V}(\rho^*,K^*)}(B_0+\penn)}{\sqrt{n\Pi_{h^*_\epsilon}}}
\Biggr)^q.
\end{eqnarray*}
With the definition of $ h^* $ and $ h^*_\epsilon$ in the beginning of
the proof, the theorem can be deduced.
\qed
\setcounter{equation}{0}
\setcounter{subsection}{0}
\begin{appendix}\label{app}

\section*{Appendix}

\subsection{Proofs of the auxiliary results}

\begin{pf*}{Proof of Proposition \protect\ref{prop full consistency}}
In this proof, we use a special case of a deviation
inequality derived in \cite{Massart07}, Corollary~6.9 (see the arXiv
version for details).
We recall that $ \firstestimator\lambda
$ is the
solution of the equation $ \tilde D_\lambda(\cdot)=0
$, thanks to the continuity of $ \rho'(\cdot) $, and we note that the
following inclusion holds:
\begin{eqnarray}
\label{eq inclusion outside} &&\bigcup_{\lambda\in\Lambda} \{\firstestimator
\lambda\notin \cF_{\delta_h^*(\lambda)} \}\nonumber\\
&&\quad \subseteq\bigcup
_{\lambda\in\Lambda} \Bigl\{\sup_{f\in\cF\backslash\cF_{\delta
_h^*(\lambda)}} \bigl\|\tilde
D_\lambda(f)-\mathbb E^0 \bigl[\tilde D_\lambda(f)
\bigr] \bigr\|_{\ell_1}\geq\inf_{f\in\cF\backslash\cF
_{\delta_h^*(\lambda)}} \bigl\|\mathbb
E^0 \bigl[\tilde D_\lambda(f) \bigr] \bigr\|_{\ell
_1} \Bigr
\}
\\
&&\quad \subseteq \Bigl\{\sup_{\lambda\in\Lambda} \Bigl[\sup_{f\in\cF
\backslash\cF_{\delta_h^*(\lambda)}}
\bigl\|\tilde D_\lambda(f)-\mathbb E^0 \bigl[\tilde
D_\lambda(f) \bigr] \bigr\|_{\ell_1}-\inf_{f\in\cF\backslash\cF
_{\delta_h^*(\lambda)}} \bigl\|
\mathbb E^0 \bigl[\tilde D_\lambda(f) \bigr]
\bigr\|_{\ell
_1} \Bigr]\geq0 \Bigr\}.\nonumber
\end{eqnarray}
Next, it holds that
%
\begin{eqnarray}
\label{eq lem decomposition criterion} \bigl\|\tilde D_\lambda(f)-\mathbb E^0 \bigl[
\tilde D_\lambda(f) \bigr] \bigr\|_{\ell_1} 
&\leq&{|\cP|} \bigl\|\tilde D_\lambda(f)-\mathbb E
\bigl[\tilde D_\lambda(f) \bigr] \bigr\|_{\ell_\infty}\nonumber\\[-8pt]\\[-8pt]
&&{}+{|\cP|} \bigl\|\mathbb E
\bigl[\tilde D_\lambda(f) \bigr]-\mathbb E^0 \bigl[\tilde
D_\lambda(f) \bigr] \bigr\|_{\ell_\infty}.\nonumber
\end{eqnarray}
By the definitions
of $\mathbb E [\tilde D_\lambda^p(\cdot) ]$ and $\mathbb
E^{0} [\tilde
D_\lambda^p(\cdot) ]$ in \eqref{def deterministic criterion}, by
change of variables and using that $\rho'(\cdot)$ is $1$-Lipschitz we
have for any
$ f\in\cF$, and any $ p\in\cP$
%
\begin{eqnarray}
\label{eq decomposition consistence} &&\sup_{f\in\cF} \bigl\|\mathbb E \bigl[\tilde
D_\lambda(f) \bigr]-\mathbb E^0 \bigl[\tilde
D_\lambda(f) \bigr] \bigr\|_{\ell_\infty}
\nonumber
\\
&&\quad  \leq\int \mu(x)K_h(x)\int \bigl\llvert \rho' \bigl(
\sigma(x)z+f^0(x)-f(x) \bigr)\nonumber\\[-8pt]\\[-8pt]
&&\hphantom{\quad  \leq\int \mu(x)K_h(x)\int \bigl\llvert} {}-\rho' \bigl(\sigma
(x)z+f^*(x)-f (x) \bigr)\bigr\rrvert \bG(z)\,\mathrm{d}z \,\mathrm{d}x
\nonumber
\\
&&\quad  \leq\me\bigl[K_h(X)\bigr] b_h.\nonumber
\end{eqnarray}
To control the stochastic term, we can then apply Massart's Inequality
to get (see the arXiv version for details)
\begin{eqnarray*}
&&\bP \biggl(\sup_{{\lambda\in\Lambda},f\in\cF}\frac{\sqrt{n\Pi
_h} \|\tilde D_\lambda(f)-\mathbb E [\tilde
D_\lambda(f) ] \|_{\ell_\infty}}{\|\rho'\|_\infty
(\sqrt{\me[\Pi_h K^2_h(X)]}+\sfrac{\|K\|_\infty}{\sqrt{n\Pi
_h}} )}\\
&&\qquad \geq 27\int
_0^{1}{H_{\cF\times\Lambda}^{1/2}(u)}\,\mathrm{d}u+{4H_{\cF}(1)}+7
\sqrt {{2z}}+2z \biggr)
\\
&&\quad \leq2|\cP|\exp(-z).
\end{eqnarray*}
Note that the factor 2 in the last inequality appears because we need to
control deviations of the absolute value of the empirical process.
Using \eqref{eq lem decomposition criterion}, \eqref{eq decomposition
consistence}, and
the last inequality, we then obtain for all $ z>0 $
%
\begin{eqnarray}
\label{eq deviations consistence} &&\bP \biggl(\sup_{\lambda\in\Lambda}\sup_{f\in\cF\backslash\cF
_{\delta_h^*(\lambda)}}
\frac{\sqrt{n\Pi_h} ( \|\tilde
D_\lambda(f)-\mathbb E^0 [\tilde
D_\lambda(f) ] \|_{\ell_1}-|\cP|{\me[K_h(X)]} b_h
)}{\|\rho'\|_\infty (\sqrt{\me[\Pi_h K^2_h(X)]}+\sfrac{\|K\|
_\infty}{\sqrt{n\Pi_h}} )}
\nonumber
\\[-8pt]\\[-8pt]
&&\quad  \geq|\cP| \biggl(27\int_0^{1}{H_{
\cF\times\Lambda}^{1/2}(u)}\,\mathrm{d}u+{4H_{\cF}(1)}+7
\sqrt{{2z}}+2z \biggr) \biggr)\leq2|\cP| e^{-z}.\nonumber
\end{eqnarray}
Now, let us have a look at $\inf_{f\in\cF\backslash\cF_{\delta
_h^*(\lambda)}}
\|\mathbb E^0 [\tilde D_\lambda(f) ] \|_{\ell_1}$ in
\eqref{eq inclusion
outside}. By the definition of $ \tilde D_\lambda(\cdot) $ and using
that ${|t^0_p-t_p|}\leq{\|t^0-t\|_{\ell_1}}$ for all $ p\in\cP$, we
have for any $
f\in\cF\backslash\cF_{\delta_h^*(\lambda)} $
\begin{eqnarray*}
\bigl\|\mathbb E^0 \bigl[\tilde D_\lambda(f) \bigr]
\bigr\|_{\ell_1} &=&\sum_{p\in\cP}\biggl\llvert \int
\biggl(\frac{x-x_0}{h} \biggr)^p\mu(x)K_h(x)\int
\rho' \bigl(\sigma(x)z+f^0(x)-f(x) \bigr) \bG(z)\,\mathrm{d}z \,\mathrm{d}x
\biggr\rrvert
\\
&\geq&\biggl\llvert \int
\frac{f^0(x)-f(x)}{\|t^0-t\|_{\ell_1}}\mu(x)K_h(x)\int \rho' \bigl(
\sigma(x)z+f^0(x)-f(x) \bigr) \bG(z)\,\mathrm{d}z \,\mathrm{d}x\biggr\rrvert ,
\end{eqnarray*}
where $\bG(\cdot)=n^{-1}\sum_{i=1}^ng_i(\cdot)$ and
$ t $ is such that $ f=\op P_t $. The last inequality is obtained using
that $ \sum_{p\in\cP}{(t^0_p-t_p)} (({x-x_0})/{h} )^p=
f(x)-f^0(x)$ and the triangular inequality. Since $
\bG(\cdot)$ is
symmetric, $ \rho'(\cdot) $ increasing (because of the convexity of $
\rho$), $
K(\cdot) $
is nonnegative, and $ \rho'(\cdot) $ is odd ($ \rho(\cdot) $ is
symmetric) and positive on $
(0,\infty) $ (because of $ \rho'(0)=0 $, the convexity of $ \rho
(\cdot) $ and the
strict convexity around 0), the last equality implies for all
$f\in\cF\backslash\cF_{\delta_h^*(\lambda)}$
\begin{eqnarray*}
&&\bigl\|\mathbb E^0 \bigl[\tilde D_\lambda(f) \bigr]
\bigr\|_{\ell_1} \\
&&\quad \geq\int \frac{ |f^0(x)-f(x) |}{\|t^0-t\|_{\ell_1}}\mu(x)K_h(x)\int
\rho' \bigl(\sigma(x)z+ \bigl|f^0(x)-f(x) \bigr| \bigr) \bG(z)\,\mathrm{d}z
\,\mathrm{d}x
\\
&&\quad \geq\int \frac{ |f^0(x)-f(x) |}{\|t^0-t\|_{\ell_1}}\mu(x)K_h(x)\int \rho'
\biggl(\sigma(x)z+\delta_h^*(\lambda)\frac{
|f^0(x)-f(x) |}{\|t^0-t\|_{\ell_1}} \biggr)
\bG(z)\,\mathrm{d}z \,\mathrm{d}x.
\end{eqnarray*}
Recall that for any $x$, $ \int\rho'(\sigma(x)z)\bG(z)\,\mathrm{d}z=0 $ thanks
to the symmetry of $ \rho(\cdot) $ and $ \bG(\cdot) $. Since $
|f^0(x)-f(x) |\|t^0-t\|_{\ell_1}^{-1}\leq1$, we obtain with
the mean value theorem for all $f\in\cF\backslash\cF_{\delta
_h^*(\lambda)}$
\begin{eqnarray*}
&&\bigl\|\mathbb E^0 \bigl[\tilde D_\lambda(f) \bigr]
\bigr\|_{\ell_1} \\
&&\quad \geq\delta_h^*(\lambda)\int \frac{ |f^0(x)-f(x) |^2}{\|t^0-t\|_{\ell_1}^2}
\mu (x)K_h(x)\inf_{u\in[0
,\delta_h^*(\lambda)] } \int \rho''
\bigl(\sigma(x)z+u \bigr) \bG(z)\,\mathrm{d}z \,\mathrm{d}x
\\
&&\quad \geq\delta_h^*(\lambda)\inf_{t:\|t\|_{\ell_1}\geq\delta
_h^*(\lambda)}\int
\frac{ |\mathrm{ P}_t(x) |^2}{\|t\|_{\ell_1}^2}\mu(x)K_h(x)\inf_{u\in[0
,\delta_h^*(\lambda)] } \int
\rho'' \bigl(\sigma(x)z+u \bigr) \bG(z)\,\mathrm{d}z \,\mathrm{d}x.
\end{eqnarray*}
We then derive, using that $
{2\delta_h^*(\lambda)}\leq\inf_{x\in V_h}\int
\rho'' (\sigma(x)z
) \bG(z)\,\mathrm{d}z$ for all $ {\lambda\in\Lambda} $ (see Condition
\ref{condition full consistency}) and $
\rho''(\cdot) $ is $\bP$-continuous,
\[
\inf_{f\in\cF\backslash\cF_{\delta_h^*(\lambda)}} \bigl\|\mathbb E^0 \bigl[\tilde
D_\lambda(f) \bigr] \bigr\|_{\ell_1} \geq\frac{\delta_h^*(\lambda)}{2}\inf
_{t:\|t\|_{\ell_1}\geq
\delta_h^*(\lambda)}\int \frac{ |\mathrm{ P}_t(x) |^2}{\|t\|_{\ell_1}^2}\mu(x)K_h(x)\int
\rho'' \bigl(\sigma(x)z \bigr) \bG(z)\,\mathrm{d}z \,\mathrm{d}x.
\]
We then observe that $ \mathrm{ P}_t(x)=t^\top U (\frac
{x-{x_0}}{h} ) $ and thus
\begin{eqnarray*}
&&\int \frac{ |\mathrm{ P}_t(x) |^2}{\|t\|_{\ell_1}^2}\mu(x)K_h(x)\int \rho''
\bigl(\sigma(x)z \bigr) \bG(z)\,\mathrm{d}z\,\mathrm{d}x
\\
&&\quad  =t^\top \biggl[\int \frac{U (\vfrac{x-{x_0}}{h} )U^\top (\vfrac
{x-{x_0}}{h} )}{\|t\|_{\ell_1}^2}\mu(x)K_h(x)\int
\rho'' \bigl(\sigma(x)z \bigr) \bG(z)\,\mathrm{d}z\,\mathrm{d}x \biggr] t.
\end{eqnarray*}
We can thus write by the definition of $\Phi_h$ in Condition \ref
{condition full consistency}
\begin{eqnarray*}
&&t^\top \biggl[\int \frac{U (\vfrac{x-{x_0}}{h} )U^\top (\vfrac
{x-{x_0}}{h} )}{\|t\|_{\ell_1}^2}\mu(x)K_h(x)\int
\rho'' \bigl(\sigma(x)z \bigr) \bG(z)\,\mathrm{d}z\,\mathrm{d}x \biggr] t\\
&&\quad \geq\frac{\|t\|_{\ell_2}^2}{\|t\|
_{\ell_1}^2}\eigen\geq \eigen/|\cP|.
\end{eqnarray*}
In summary, we have for any $\lambda\in\Lambda$,
$
\inf_{f\in\cF\backslash\cF_{\delta_h^*(\lambda)}} \|\mathbb
E^0 [\tilde
D_\lambda(f) ] \|_{\ell_1}
\geq\frac{\eigen\delta_h^*(\lambda)}{2{|\cP|}}.
$
By the definition of $ \delta_h^*(\lambda) $ in Condition \ref
{condition full consistency} and as $ n\Pi_h\geq1 $, it holds that
\[
\delta_h^*(\lambda)>2|\cP|^2\frac{\|\rho'\|_\infty (\sqrt {\me[\Pi_h K^2_h(X)]}+\sfrac{\|K\|_\infty}{\sqrt{n\Pi_h}}
)}{\eigen\sqrt{n\Pi_h} (E^*+7\sqrt{{4\ln(2|\cP|n)}}+4\ln
(2|\cP|n) )^{-1}}+2|
\cP|^2{\me\bigl[K_h(X)\bigr]} \frac
{b_h}{\eigen}.
\]
Using Inequalities \eqref{eq inclusion
outside} and
\eqref{eq deviations
consistence} with $ z=\ln(2|\cP|n) $, and the last inequality, we obtain
\begin{eqnarray*}
\bP \biggl(\bigcup_{\lambda\in\Lambda} \{\firstestimator \lambda
\notin\cF_{\delta_h^*(\lambda)} \} \biggr) &\leq& \bP \biggl(\sup_{\lambda\in\Lambda}
\sup_{f\in\cF\backslash\cF
_{\delta_h^*(\lambda)}} \biggl[ \bigl\|\tilde D_\lambda(f)-\mathbb
E^0 \bigl[\tilde D_\lambda(f) \bigr] \bigr\|_{\ell_1}-
\frac{\eigen\delta_h^*(\lambda
)}{2{|\cP|}} \biggr]\geq0 \biggr)\\
 &\leq&1/n^2.
\end{eqnarray*}
\upqed
\end{pf*}
\begin{pf*}{Proof of Proposition \protect\ref{prop uniform deviations}}
The definitions of
$\hat f_\lambda$ and $f^0$ (see (\ref{def estimator with
fixed lambda}) and (\ref{coefficient_taylor}), resp.) imply that
$ |\hat{f}_\lambda({x_0})-f^*({x_0}
) |= | (\hat t_\lambda )_ {
0,\ldots,0
} -t^0_{0, \ldots,0} |\nonumber
\leq\|\hat t_\lambda-t^0 \|_{\ell_\infty}.
$
Using $
\hat f_\lambda\in\cF_{\delta_h^*(\lambda)} $, Lemma~\ref{lem
criterion inversible}, and the
last inequality, we have
\begin{eqnarray*}
\bigl|\hat{f}_\lambda({x_0})-f^*({x_0} ) \bigr|\leq
\tfrac{4}{3}c_\lambda^{-1} \bigl\|\mathbb E^0
\bigl[\tilde D_\lambda(\hat{f}_\lambda) \bigr]-\mathbb
E^0 \bigl[\tilde D_\lambda\bigl(f^0\bigr) \bigr]
\bigr\|_{\ell_\infty}.
\end{eqnarray*}
Recall that by definition $ \tilde D_\lambda(\hat{f}_\lambda)=0 $
and $
\mathbb E^0 [\tilde
D_\lambda(f^0) ]=0 $. Thus, for all $\lambda\in\Lambda$ such
that $ \hat
f_\lambda\in
\cF_{\delta_h^*(\lambda)} $, the last inequality implies
\begin{eqnarray*}
\bigl|\hat{f}_\lambda({x_0})-f^*({x_0} ) \bigr| &\leq&
\tfrac{4}{3}c_\lambda^{-1} \bigl(\bigl\|\tilde
D_\lambda(\hat{f}_\lambda)-\mathbb E \bigl[\tilde
D_\lambda(\hat{f}_\lambda) \bigr]\bigr\|_{\ell_\infty}+\bigl\|\mathbb E
\bigl[\tilde D_\lambda(\hat{f}_\lambda) \bigr]-\mathbb
E^0 \bigl[\tilde D_\lambda(\hat{f}_\lambda) \bigr]
\bigr\|_{\ell_\infty} \bigr).
\end{eqnarray*}
From Lemma~\ref{lem control bias} and the last display, we obtain
\begin{eqnarray*}
\bigl|\hat{f}_\lambda({x_0})-f^*({x_0} ) \bigr| &\leq&
\tfrac{4}{3}c_\lambda^{-1} \bigl(\bigl\|\tilde
D_\lambda(\hat{f}_\lambda)-\mathbb E \bigl[\tilde
D_\lambda(\hat{f}_\lambda) \bigr]\bigr\|_{\ell_\infty}+
\tfrac{5}{4}{c} _\lambda b_h \bigr)
\\
&\leq&\tfrac{5}{3} b_h+\tfrac{4}{3}\sup
_{\lambda\in\Lambda}\sup_{f\in\cF_{\delta_h^*(\lambda)}}c_\lambda^{-1}
\bigl\|\tilde D_\lambda(f)-\mathbb E \bigl[\tilde D_\lambda(f) \bigr]
\bigr\|_{\ell_\infty}.
\end{eqnarray*}
This yields
\[
\bigl|\hat{f}_\lambda({x_0})-f^*({x_0} ) \bigr|
\leq3b_h+2\sup_{\lambda\in\Lambda}\sup_{f\in\cF_{\delta
_h^*(\lambda)}}c_\lambda^{-1}
\bigl\|\tilde D_\lambda(f)-\mathbb E \bigl[\tilde D_\lambda(f) \bigr]
\bigr\|_{\ell_\infty}.
\]
From the last inequality and the definitions of $ \mathrm{V}(\cdot) $ and $
c_\lambda$ introduced in \eqref{def Huber Variance} and \eqref{def
normalization term}, respectively, we deduce
\begin{eqnarray*}
&&\bP \biggl( \biggl\{\sup_{\lambda\in\Lambda} \biggl[ \bigl|
\hat{f}_{
\lambda}({x_0} )-f^*(x_0) \bigr|-2
\frac{\sqrt{\mathrm{V}(\lambda)} B_z}{\sqrt{n\Pi
_h}} \biggr]\geq 3b_h \biggr\}\cap\bigcap
_{\lambda\in\Lambda} \{\hat f_\lambda \in \cF_{\delta_h^*(\lambda)} \}
\biggr)
\\
&&\quad \leq\bP \biggl(\sup_{\lambda\in\Lambda}\sup_{f\in\cF_{\delta
_h^*(\lambda)}}
\biggl[2c_\lambda^{-1} \bigl\|\tilde D_\lambda(f)-\mathbb E
\bigl[\tilde D_\lambda(f) \bigr]\bigr\|_{\ell_\infty}-2\frac{\sqrt{\mathrm{V}(\lambda)}
B_z}{\sqrt{
n\Pi_h}}
\biggr]\geq0 \biggr)
\\
&&\quad \leq\bP \biggl(\sup_{\lambda\in\Lambda}\sup_{f\in\cF_{\delta
_h^*(\lambda)}}
\frac{
\|\tilde
D_\lambda(f)-\mathbb E [\tilde D_\lambda(f) ]\|_{
\ell_\infty}}{\sqrt{\Pi_h}\sqrt{\mathbb EP_n [\lambda
'(f^*) ]^2}+
\lambda'_\infty\sfrac{\ln^2(n)}{\sqrt{n\Pi_h}}}\geq \frac{B_z}{\sqrt{n\Pi_h}} \biggr).
\end{eqnarray*}
Using Proposition~\ref{prop large deviation} and the last inequality,
we finally obtain
\begin{eqnarray*}
\label{eq prop control first event} \bP \biggl( \biggl\{\sup_{\lambda\in\Lambda} \biggl[ \bigl|
\hat{f}_{
\lambda}({x_0} )-f^*(x_0) \bigr|-2
\frac{\sqrt{\mathrm{V}(\lambda)} B_z}{\sqrt{n\Pi
_h}} \biggr]\geq 3b_h \biggr\}\cap\bigcap
_{\lambda\in\Lambda} \{\hat f_\lambda \in \cF_{\delta_h^*(\lambda)} \}
\biggr)\leq 2|\cP|\mathrm{e}^{-z}.
\end{eqnarray*}
\upqed
\end{pf*}
\begin{pf*}{Proof of Proposition \protect\ref{prop control huber variance}}
We first recall by the definition of the estimator \eqref{def adaptive lambda}
\[
\sqrt{\widehat{\mathrm{V}} (\lambda)}=\frac{\sqrt {\Pi_h P_n [\lambda'(\hat f_{\lambda}) ]^2}+
\lambda'_\infty\sfrac{\ln^2(n)}{\sqrt{n\Pi_h}}}{P_n \lambda
''(\hat f_{\lambda})},
\]
where
\[
\Pi_h P_n \bigl[\lambda'(\hat
f_{\lambda}) \bigr]^2=\sum_{i=1}^n
\frac{1}{n\Pi_h} \bigl[\rho '\bigl(Y_i-\hat
f_{\lambda}(X_i)\bigr) \bigr]^2K^2
\biggl(\frac{X_i-x_0}{h} \biggr)
\]
and
\[
P_n \lambda''(\hat f_{\lambda})=
\sum_{i=1}^n\frac{1}{n\Pi_h}\rho
''\bigl(Y_i-\hat f_{\lambda}(X_i)
\bigr)K \biggl(\frac{X_i-x_0}{h} \biggr).
\]
In the following, we assume to be on the event
$\bigcap_{\lambda\in\Lambda} \{\hat f_\lambda\in
\cF_{\delta_h^*(\lambda)} \}$, which is true with probability
at least $1-1/n^2$ according to Proposition~\ref{prop full
consistency}. Then, using Massart's Inequality (see the arXiv version
for details)
we can control the deviation of the process $\Pi_h P_n [\lambda
'(\hat
f_{\lambda}) ]^2$ as follows:
%
\begin{equation}
\label{eq lem control variance term} \bP \biggl(\sup_{\lambda\in\Lambda}\sup_{f\in\cF_{\delta
_h^*(\lambda)}}
\bigl(\lambda'_\infty\bigr)^{-2} \bigl\llvert
\Pi_hP_n \bigl[\lambda'(f)
\bigr]^2-\Pi_h\mathbb EP_n \bigl[
\lambda'(f) \bigr]^2\bigr\rrvert \geq\frac{B_{2\ln(n)}}{\sqrt{n\Pi_h}}
\biggr)\leq2/n^2,
\end{equation}
where $B_{\cdot}$ is defined in \eqref{def upper bound Massart}.
Similarly, using again Massart's Inequality,
we control the deviation of $P_n \lambda''(\hat f_{\lambda})$ as follows:
%
\begin{equation}
\label{eq lem control denominator} \bP \biggl(\sup_{\lambda\in\Lambda}\sup_{f\in\cF_{\delta
_h^*(\lambda)}}
\|K\|_\infty^{-1}\bigl\llvert P_n
\lambda''(f)-\mathbb EP_n\lambda
''(f)\bigr\rrvert \geq\frac{B_{2\ln(n)}}{\sqrt{n\Pi_h}} \biggr)
\leq2/n^2.
\end{equation}
Then, for any $ \lambda\in\Lambda$, by the continuity of $ \rho' $
and $ \rho'' $ almost everywhere,
$\|\rho''\|_\infty\leq1$, and the mean
value theorem,
we have for all $ f\in\cF_{\delta_h^*(\lambda)} $
\begin{eqnarray*}
&&\Pi_h \bigl\llvert \mathbb EP_n \bigl[
\lambda'(f) \bigr]^2-\mathbb EP_n \bigl[
\lambda'\bigl(f^*\bigr) \bigr]^2\bigr\rrvert \\
&&\quad \leq
\frac{1}{n\Pi_n}\sum_{i=1}^n\mathbb{E}
\bigl\llvert \rho'\bigl(Y_i-f(X_i)
\bigr)^2-\rho'\bigl(Y_i-f^*(X_i)
\bigr)^2\bigr\rrvert K^2 \biggl(\frac{X_i-x_0}{h}
\biggr)
\\
&&\quad \leq2\|K\|_\infty^2\bigl(\delta_h^*(
\lambda)+b_h\bigr).
\end{eqnarray*}
Similarly,
$ \sup_{f\in\cF_{\delta_h^*(\lambda)}}
\llvert \mathbb EP_n\lambda''(f)-\mathbb EP_n\lambda''(f^*)\rrvert \leq
\|K\|_\infty(\delta_h^*(\lambda)+b_h).
$
Note that for any $ \lambda\in\Lambda$
\[
s_n:=\|K\|_\infty s_h(\lambda)\geq \bigl(1\vee2\|K
\|_\infty \bigr)\|K\|_\infty \bigl[\delta_h^*(
\lambda)+b_h \bigr]+ \bigl[\bigl(\lambda '_\infty
\bigr)^{2}\vee\|K\|_\infty \bigr]\frac{B_{2\ln(n)}}{\sqrt {n\Pi_h}},
\]
and we observe (under Condition \ref{condition variance estimation})
that $
s_n\leq\frac{1}{2}\min \{\mathbb EP_n\lambda''(f^*),\Pi_h
\mathbb EP_n [\lambda'(f^*) ]^2 \}.
$
Using this, \eqref{eq lem control variance term}, and \eqref{eq lem
control denominator}, we obtain with probability $ 1-5/n^2 $ for any $
\lambda\in\Lambda$
\[
\sqrt{\widehat{\mathrm{V}} (\lambda)}\leq\frac{\sqrt {\Pi_h \mathbb EP_n [\lambda'(f^*) ]^2+s_n}+
\lambda'_\infty\sfrac{\ln^2(n)}{\sqrt{n\Pi_h}}}{\mathbb
EP_n\lambda''(f^*)-s_n}\leq\sqrt6\sqrt{{\mathrm{V}}
(\lambda)}
\]
and
\[
\sqrt{\widehat{\mathrm{V}} (\lambda)}\geq\frac{\sqrt {\Pi_h \mathbb EP_n [\lambda'(f^*) ]^2-s_n}+
\lambda'_\infty\sfrac{\ln^2(n)}{\sqrt{n\Pi_h}}}{\mathbb
EP_n\lambda''(f^*)+s_n}\geq\frac{\sqrt2}{3}
\sqrt{{\mathrm{V}} (\lambda)}.
\]
(Instead of the given factors in front of $\sqrt{{\mathrm{V}}
(\lambda)}$, one could readily obtain factors that tend to one as
$n\to\infty$. This is of minor interest here.) This proves the
claim.
\end{pf*}

\subsection{Technical lemmas}

We first give a result for the {
deterministic criterion} $ \mathbb E^0 [\tilde D_\lambda(\cdot
) ] $ defined in
(\ref{def deterministic criterion}):
%
\begin{lemma}\label{lem criterion inversible}
Let $ \lambda$ be as in \eqref{def estimator with fixed lambda},
$n\in\{1,2,\dots\}$, and
$h\in(0,1]^d$ such that Condition \ref{condition full consistency} is
satisfied,
the following holds:\vadjust{\goodbreak}
\begin{enumerate}[2.]
\item$ \mathbb E^0 [\tilde
D_\lambda(f^{0}) ]=0 $,
and the function $ \mathbb E^0 [\tilde D_\lambda(f) ] $ is
bijective as
function of
$
\cF_{\delta_h^*(\lambda)} $ (see Definition \eqref{def event
consistence}) on the corresponding image.
\item For any $ f,\tilde
f\in\cF_{\delta_h^*(\lambda)} $,
$
\|t-\tilde t\|_{\ell_\infty}
\leq
\tfrac{4}{3}c_\lambda^{-1}
\|\mathbb E^0 [\tilde D_\lambda(f) ]-\mathbb E^0 [\tilde
D_\lambda(\tilde f) ]\|_{\ell_\infty},
$
where $ \mathrm{P}_t=f $ and $ \mathrm{P}_{\tilde t}=\tilde f $.
\end{enumerate}
\end{lemma}

Next, we consider the bias.
%
\begin{lemma}\label{lem control bias}
Let $ \lambda$ be as in \eqref{def estimator with fixed lambda},
$n\in\{1,2,\dots\}$, and
$h\in(0,1]^d$ such that Condition \ref{condition full consistency} is
satisfied, it holds that
\[
\sup_{f\in\cF_{\delta_h^*(\lambda)}} \bigl\|\mathbb E^0 \bigl[\tilde
D_\lambda(f) \bigr]-\mathbb E \bigl[\tilde D_\lambda(f) \bigr]
\bigr\|_{\ell_\infty}\leq\tfrac{5}{4}{c}_\lambda b_h.
\]
\end{lemma}

Next, we do some simple algebra.
%
\begin{lemma}\label{lemma decomposition}
For any $x,y\in[0,\infty)$, it holds that
$
x^q\leq2^q[x-y]_+^q+2^qy^q.
$
Moreover, for any $l,q\in\{1,2,\dots\}$ and $x_1,\ldots,x_l\geq0$,
it holds,
that
$
(\sum_{i=1}^lx_i )^q\leq l^{q-1} (\sum_{i=1}^lx_i^q ).
$
\end{lemma}
The proof consists of simple algebra and is available in the arXiv version.

The following lemma allows us to
get our hands on the estimator $ \widehat{\mathrm{V}}(\cdot) $.
%
\begin{lemma}\label{lem aniso criterion inversible} Let $ \cD_h(\cdot
)\dvtx [-M,M]\rightarrow\bR$ and $ c_{\hat\rho} $ be as defined in
the proof of Theorem~\ref{Th adaptive aniso} and assume $ f^*\in
\mathbb{H}_d(\vec{\beta},L,M) $ and $n$ sufficiently large such that
Condition \ref{condition full consistency} is satisfied for all $ h\in
\cH$. Then, for any
$h\in\cH$ and $ t,\tilde t\in[f^*(x_0)-\delta_h^*(\lambda
),f^*(x_0)+\delta_h^*(\lambda)] $, it holds that $
|t-\tilde t|\leq\frac{4}{3}c_{\hat\rho}^{-1}|\cD_h(t)-\cD
_h(\tilde t)|.
$
\end{lemma}

The proof of the lemma is similar to the one of Lemma~\ref{lem
criterion inversible} (see the arXiv version for details).

Next, we control the distance of $\tilde\cD_{h}(f)$ to $\cD_{h}(f)$ for
appropriate bandwidth $h$ and functions~$f$:
%
\begin{lemma}\label{lem deviation criterion anisotrope} For n
sufficiently large $n$ sufficiently large such that Conditions \ref
{condition full consistency} and \ref{condition variance estimation}
are satisfied for all $ h\in\cH$. It holds that
\[
\mathbb Ec_{\hat\rho}^{-q}\sup_{h\in\cH: h\succeq
h^*_\epsilon}\sup
_{f\in\cF_{\tilde\delta_h}}\bigl\llvert \tilde\cD_{h}(f)-
\cD_{h}(f)\bigr\rrvert ^q \asymp2^q \biggl(
\frac{\sqrt{6\mathrm{V}(\rho^*,K^*)}(B_0+\penn
)}{\sqrt{n\Pi_{h_\epsilon^*}}} \biggr)^q,
\]
where $ \tilde\delta_h $, $ h_\epsilon^* $, $ \tilde\cD$ and $\cD
$ are defined in the proof of Theorem~\ref{Th adaptive aniso},
$\operatorname{Gamma}(q)$ is the classical Gamma function, $ \mathrm{V}(\rho^*,K^*) $ is
defined in \eqref{def homo variance} and \eqref{def new oracle
aniso}, $ \penn$ is defined in Section~\ref{section full adaptive
anisotropic}.
\end{lemma}
The proof is an application of Proposition~\ref{prop large deviation}
(see the arXiv version for details).

Eventually, we look at the distance to $\cD_{{h'}\vee{h}}(f)$ to
$\cD_{{h}}(f)$ for
appropriate bandwidths $h$ and $h'$ and functions $f$:
%
\begin{lemma}\label{lem control bias anisotrope}
For any $ f^*\in
\mathbb{H}_d(\vec{\beta},L,M) $ such that $ \vec\beta\in(0,1]^d
$, and $n$ sufficiently large such that Condition \ref{condition full
consistency} is satisfied for all $ h\in\cH$,
it holds that for any $ {h, h'}\in\cH$
\[
\sup_{f\in\cF_{\tilde\delta_h}} \bigl|\cD_{{h'}\vee{
h}}(f)-\cD_{{h}}(f) \bigr|
\leq\frac{5} 4c_{\hat\rho}L\sum_{j=1}^d
\bigl(h_j'\bigr)^{\beta_j},
\]
where $ \cD_h $ and $ c_{\hat\rho} $ are defined in \eqref{def re
partial expectation anisotrope} and \eqref{def re denominator} in the
proof of Theorem~\ref{Th adaptive aniso}.
\end{lemma}

\subsection{Proofs of the technical lemmas}

\begin{pf*}{Proof of Lemma \protect\ref{lem criterion inversible}}
Let us proof the first claim. For this, we note that the components of $
\mathbb E^0 [\tilde
D_\lambda(f) ]$ are given by
\begin{eqnarray*}
\mathbb E^0 \bigl[\tilde D_\lambda^p(f) \bigr]=
\int \biggl(\frac{x-x_0}{h} \biggr)^p\mu(x)K_h(x)\int
\rho' \bigl(\sigma(x)z+f^0(x)-f(x) \bigr)
\frac{1}{n}\sum_{i=1}^ng_i(z)\,\mathrm{d}z
\,\mathrm{d}x.
\end{eqnarray*}
Since $ \rho(\cdot) $ and $\sum_i g_i(\cdot)$ are
symmetric, it holds that $ \int
\rho'(z)\sum_i g_i(z)\,\mathrm{d}z=0 $ and $ \mathbb E^0 [\tilde D_\lambda
^p(f^0) ]=0
$. We now show that $ \mathbb E^0 [\tilde D_\lambda^p(\cdot)
] $ is injective
on the image of
$ \cF_{\delta_h^*(\lambda)} $ exploiting further the symmetry of $
\rho(\cdot) $ and
$\sum_i g_i(\cdot)$. Consider $ f,\tilde f\in\cF_{\delta
_h^*(\lambda)}$
such
that $ \mathbb E^0 [\tilde D_\lambda(f) ]=\mathbb E^0
[\tilde
D_\lambda(\tilde f) ] $. We have to show that $
f=\tilde f $. For this, we first note that
\[
\sum_{p\in\cP}(t_p-\tilde t_p)
\bigl(\mathbb E^0 \bigl[\tilde D_\lambda^p(
\mathrm{P}_t) \bigr]-\mathbb E^0 \bigl[\tilde
D_\lambda^p(\mathrm{P}_{\tilde t}) \bigr] \bigr)=0,
\]
where $ t $ and $ \tilde t $ are such that $ \mathrm{ P}_t=f $ and $ \mathrm{
P}_{\tilde
t}= \tilde f$.
To simplify the presentation, we introduce the notation $
u(\cdot):=(f-f^0)(\cdot)$, $
\tilde u(\cdot):=(\tilde f-f^0)(\cdot)$, and $
\bG(\cdot):=n^{-1}\sum_{i=1}^ng_i(\cdot) $. Since $
\bG(\cdot)$ is
symmetric, $ K(\cdot) $ is nonnegative, and $ \rho'(\cdot) $ is odd
and positive on $
(0,\infty)$, the last display implies
\begin{eqnarray*}
&&\int K_h(x)\mu(x) \bigl[u(x)-\tilde u(x) \bigr]\\
&&\quad {}\times\int \bigl[
\rho' \bigl(\sigma(x)z-u(x) \bigr)-\rho' \bigl(\sigma
(x)z-\tilde u(x) \bigr) \bigr] \bG(z)\, \mathrm{d}z \,\mathrm{d}x=0
\\
&&\qquad \Leftrightarrow\qquad \int K_h(x)\mu(x) \bigl|u(x)-\tilde u(x) \bigr|\\
&&\hphantom{\qquad \Leftrightarrow\qquad}\quad  {}\times\int \bigl
\llvert \rho' \bigl(\sigma(x)z-u(x) \bigr)-\rho'
\bigl(\sigma (x)z-\tilde u(x) \bigr)\bigr\rrvert \bG(z) \,\mathrm{d}z \,\mathrm{d}x=0.
\end{eqnarray*}
As $ f,\tilde f\in\cF_{\delta_h^*(\lambda)} $, it holds that $
\sup_{x\in
V_h}|u(x)|\vee|\tilde u(x)|\leq\delta_h^*(\lambda) $. Moreover,
using the mean value
theorem,
the $\bP$-continuity of $ \rho'' $ and Condition \ref{condition full
consistency}, we obtain
\begin{eqnarray*}
&&\int K_h(x)\mu(x) \bigl|u(x)-\tilde u(x) \bigr|\int \bigl\llvert
\rho' \bigl(\sigma(x)z-u(x) \bigr)-\rho' \bigl(\sigma
(x)z-\tilde u(x) \bigr)\bigr\rrvert \bG(z) \,\mathrm{d}z \,\mathrm{d}x
\\
&&\quad \geq\int K_h(x)\mu(x) \bigl|u(x)-\tilde u(x) \bigr|^2\inf
_{s:|s|\leq\delta_h^*(\lambda)}\int \rho'' \bigl(
\sigma(x)z-s\bigr) \bG(z)\,\mathrm{d}z \,\mathrm{d}x
\\
&&\quad \geq\int K_h(x)\mu(x) \bigl|u(x)-\tilde u(x) \bigr|^2\inf
_{s:|s|\leq\delta_h^*(\lambda)}\int \rho'' \bigl(
\sigma(x)z-s\bigr) \bG(z)\,\mathrm{d}z \,\mathrm{d}x
\\
&&\quad \geq\int K_h(x)\mu(x) \bigl|u(x)-\tilde u(x) \bigr|^2 \biggl[
\int \rho'' \bigl(\sigma(x)z\bigr) \bG(z)\,\mathrm{d}z-
\delta_h^*(\lambda) \biggr]\,\mathrm{d}x
\\
&&\quad \geq\frac{1}{2}\int K_h(x)\mu(x) \bigl|u(x)-\tilde u(x)
\bigr|^2\int\rho'' \bigl(\sigma(x)z\bigr)
\bG(z)\,\mathrm{d}z \,\mathrm{d}x.
\end{eqnarray*}
The last display, Condition \ref{condition full consistency}, and the
nonnegativity of $ K(\cdot) $ over its support yield that there exists
an nonempty open set $ \cV$ such that
$
\sup_{x\in\cV} |u(x)-\tilde u(x) |=0.
$
As $ u $ and $ \tilde u $ are polynomials with finite degree, we
finally obtain
that $ f=\tilde f $, and the first claim is proved.

Let us now turn to the second claim. We set $ D(\cdot):=\mathbb
E^0 [\tilde
D_\lambda(\cdot) ] $ and note that $ D(\cdot)$ is differentiable and
injective on $ \cF_{\delta_h^*(\lambda)}$ (the latter according to
the first claim). We can
consequently find an inverse of the function $D(\cdot)$ on the image
of
$D(\cdot)$ on $ \cF_{\delta_h^*(\lambda)} $. We then obtain,
denoting the\vspace*{-1pt} matrix $
\ell_\infty$-norm by $ |\hspace*{-1.1pt}|\hspace*{-1.1pt}|\cdot|\hspace*{-1.1pt}|\hspace*{-1.1pt}|_\infty$ and the inverse of
$D(\cdot)$ by
$D^{-1}(\cdot)$, for all $ f\in\cF_{\delta_h^*(\lambda)} $
\[
\bigl|\hspace*{-1.1pt}\bigl|\hspace*{-1.1pt}\bigl|J_{D^{-1}}(f)\bigr|\hspace*{-1.1pt}\bigr|\hspace*{-1.1pt}\bigr|_\infty=
\bigl|\hspace*{-1.1pt}\bigl|\hspace*{-1.1pt}\bigl|J_D^{-1}
(f)\bigr|\hspace*{-1.1pt}\bigr|\hspace*{-1.1pt}\bigr|_\infty=\bigl|\hspace*{-1.1pt}\bigl|\hspace*{-1.1pt}\bigl|J_D(f)
\bigr|\hspace*{-1.1pt}\bigr|\hspace*{-1.1pt}\bigr|_\infty^{-1}
\leq \bigl[ J_D(f) \bigr]_{0,0}^{-1}= \bigl[
\mathbb EP_n\lambda''(f)
\bigr]^{-1} \leq\tfrac{4}{3}c_{\lambda}^{-1}.
\]
The constant $ c_{\lambda} $ is defined in \eqref{def normalization
term} and the last inequality is obtained by the $\bP$-continuity of $
\rho''(\cdot)$
and Condition \ref{condition full consistency}.\vspace*{-1pt}
The mean value theorem and the last inequality then imply
for any $ f,\tilde f\in\cF_{\delta_h^*(\lambda)}$ and the associated
coefficients $ t $ and $ \tilde t $
\[
\|t-\tilde t\|_{\ell_\infty}=\bigl\llVert D^{-1}\circ
D(f)-D^{-1}\circ D(\tilde f)\bigr\rrVert _{\ell_\infty}\leq
\tfrac{4}{3}c_{\lambda}^{-1}\bigl\llVert D(f)-D(\tilde f)
\bigr\rrVert _{\ell_\infty}.
\]
This proves the second claim.
\end{pf*}
%
\begin{pf*}{Proof of Lemma \protect\ref{lem control bias}}
By the definitions
of $\mathbb E [\tilde D_\lambda^p(\cdot) ]$ and $\mathbb
E^0 [\tilde
D_\lambda^p(\cdot) ]$ in \eqref{def deterministic criterion}, we
have for any
$ f\in\cF_{\delta_h^*(\lambda)} $, any $ \lambda\in\Lambda$ ,
and any $ p\in\cP$
%
\begin{eqnarray}
\label{eq lem first decomposition}
 && \bigl|\mathbb E^0 \bigl[\tilde D_\lambda^p(f)
\bigr]-\mathbb E \bigl[\tilde D_\lambda^p(f) \bigr] \bigr|
\nonumber
\\
&&\quad  \leq\int \mu(x)K_h(x)\\
&&\qquad {}\times\int \bigl\llvert \rho' \bigl(
\sigma(x)z+f^0(x)-f(x) \bigr)-\rho' \bigl(\sigma
(x)z+f^*(x)-f (x) \bigr)\bigr\rrvert \bG(z)\,\mathrm{d}z \,\mathrm{d}x.\nonumber
\end{eqnarray}
It additionally holds for all $f\in\cF_{\delta_h^*(\lambda)}$ that
$ \sup_{x\in V_h
}|f^0(x)-f(x)|\leq\delta_h^*(\lambda) $. Together with the
definition of $ f^0 $ in
\eqref{coefficient_taylor}, this implies for any $f\in\cF_{\delta
_h^*(\lambda)} $
\[
\sup_{x\in V_h
}\bigl|f^*(x)-f(x)\bigr|\leq\sup_{x\in V_h
}\bigl|f^*(x)-f^0(x)\bigr|+
\sup_{x\in V_h
}\bigl|f^0(x)-f(x)\bigr|\leq b_h+
\delta_h^*(\lambda).
\]
This implies, due to the mean value theorem, that there is a
$u_x\in\bR: |u_x|\leq b_h+\delta_h^*(\lambda) $ such that
\begin{eqnarray*}
\bigl\llvert \rho' \bigl(\sigma(x)z+f^0(x)-f(x)
\bigr)-\rho' \bigl(\sigma (x)z+f^*(x)-f (x) \bigr)\bigr\rrvert
\leq\bigl|f^*(x)-f^0(x)\bigr|\rho'' \bigl(
\sigma(x)z+u_x \bigr).
\end{eqnarray*}
Using Condition \ref{condition full consistency}, \eqref{eq
lem first
decomposition}, the last inequality, and the
definitions $ b_h $, and $ c_\lambda$ defined in \eqref{def bias
term} and \eqref{def
normalization term} respectively, we obtain for any $ \lambda\in
\Lambda$
\begin{eqnarray*}
&&\sup_{f\in\cF_{\delta_h^*(\lambda)}} \bigl\|\mathbb E^0 \bigl[\tilde
D_\lambda(f) \bigr]-\mathbb E \bigl[\tilde D_\lambda(f) \bigr]
\bigr\|_{\ell_\infty}
\\
&&\quad \leq \int\mu(x)K_h(x)\bigl|f^*(x)-f^0(x)\bigr|\int \bigl[
\rho'' \bigl(\sigma(x)z \bigr)+b_h+
\delta_h^*(\lambda) \bigr] \bG(z)\,\mathrm{d}z \,\mathrm{d}x\leq\frac{5}{4}{c}_\lambda
b_h.
\end{eqnarray*}
\upqed
\end{pf*}
\begin{pf*}{Proof of Lemma \protect\ref{lem control bias anisotrope}}
Recall that we consider the uniform design and the homoscedastic noise level.
By the definition of $ \cD_{h} $ and with a change of variables, we have
\begin{eqnarray*}
&&\sup_{f\in\cF_{\tilde\delta_h}} \bigl|\cD_{{h'}\vee{
h}}(f)-\cD_{{h}}(f) \bigr|\\
 &&\quad =
\sup_{f\in\cF_{\tilde\delta_h}}\biggl\llvert \int \hat K(x)\int \hat
\rho' \bigl(\sigma z+f^*\bigl(x_0+{h}\vee
{h}'x\bigr)-f(x_0) \bigr)g(z) \,\mathrm{d}z \,\mathrm{d}x
\\
&&\hphantom{\quad =
\sup_{f\in\cF_{\tilde\delta_h}}\biggl\llvert}{} -\int \hat K(x)\int\hat\rho' \bigl(\sigma
z+f^*(x_0+{h}x)-f(x_0) \bigr)g(z) \,\mathrm{d}z \,\mathrm{d}x\biggr\rrvert .
\end{eqnarray*}
Using $ f\in\bH_d(\vec\beta,L,M) $, the $\bP$-continuity of $ \rho
''(\cdot) $, the last equality, and the mean value theorem, we obtain:
\begin{eqnarray*}
&&\sup_{f\in\cF_{\tilde\delta_h}} \bigl|\cD_{{h'}\vee{h}}(f)-\cD _{{h}}(f) \bigr| \\
&&\quad
\leq\sup_{|s|\leq\tilde\delta_h+b_{h}}\int\hat\rho''(\sigma
z+s)g(z)\,\mathrm{d}z\int\hat K(x) \bigl|f^*\bigl(x_0+{h}\vee{h}'x
\bigr)-f^*(x_0+{h}x) \bigr|\,\mathrm{d}x
\\
&&\quad \leq \biggl(\int\hat\rho''(\sigma z)g(z)\,\mathrm{d}z+\tilde
\delta _h+b_{h} \biggr)L\sum_{j=1}^d\bigl|{h_j}
\vee{h_j'}-h_j\bigr|^{\beta_j}.
\end{eqnarray*}
With Condition \ref{condition full consistency} and definition of $
\tilde\delta_h $ in Proof of Theorem~\ref{Th adaptive aniso}, this yields
\[
\sup_{f\in\cF_{\tilde\delta_h}} \bigl|\cD_{{h'}\vee{h}}(f)-\cD _{{h}}(f) \bigr|
\leq\frac{5} 4 c_{\hat\rho}L\sum_{j=1}^d
\bigl(h_j'\bigr)^{\beta_j}.
\]
\upqed
\end{pf*}
\end{appendix}

\section*{Acknowledgements}
The authors acknowledge partial financial support as members of the
German--Swiss Research Group FOR916 (Statistical
Regularization and Qualitative Constraints) with grant number 20PA20E-134495/1.

We thank Oleg Lepski, Joseph Salmon, and Sara van de Geer for many
helpful discussions.
We also thank the associate editor and the referees for their valuable comments.



\printhistory

\end{document}